\documentstyle[12pt]{article}

\font\twelvefrak=eufm10 at 12pt

\font\sevenfrak=eufm7
\font\fivefrak=eufm5
\newfam\frakfam
\textfont\frakfam=\twelvefrak
\scriptfont\frakfam=\sevenfrak
\scriptscriptfont\frakfam=\fivefrak
\def\frak{\fam\frakfam\twelvefrak}

\font\twelveBbb=msbm10 at 12pt

\font\sevenBbb=msbm7
\font\fiveBbb=msbm5
\newfam\Bbbfam
\textfont\Bbbfam=\twelveBbb
\scriptfont\Bbbfam=\sevenBbb
\scriptscriptfont\Bbbfam=\fiveBbb
\def\Bbb{\fam\Bbbfam\twelveBbb}

\newtheorem{lem}{Lemma.}

\newtheorem{prop}{Proposition.}

\newtheorem{cor}{Corollary.}

\newtheorem{th}{Theorem.}

\newtheorem{thintro}{Theorem}

\newtheorem{numth}{Theorem}[subsection] 

\newtheorem{propC}{Proposition (C).}

\newtheorem{propAR}{Proposition (AR).}

\def\pf{{\it Proof. }}

\def\adots
{\mathinner{\mkern2mu\raise1pt\hbox{.}\mkern3mu\raise4pt\hbox{.}\mkern1mu\raise7pt\hbox{.}}}

\def\lra{\longrightarrow} 

\def\binom#1#2{{\left({#1 \atop #2}\right)}}

\def\dist{{\rm Dist}}
\def\Ext{\mbox{\rm Ext}} 
\def\gr{\mathop{\rm gr}\nolimits}
\def\grisol{{\rm gr}_{\rm isol}}
\def\Hom{\mbox{Hom}}
\def\Ind{\mbox{Ind}}
\def\Ker{\mathop{\rm Ker}}
\def\Ima{\mathop{\rm Im}}
\def\Lie{\mathop{\rm Lie}}
\def\Max{\mbox{Max}} 
\def\Tor{\mbox{\rm Tor}}

\def\frakg{{\frak g}}
\def\frakgl{{\frak gl}}
\def\frakm{{\frak m}}

\def\frakp{{\frak p}}

\def\frakt{{\frak t}}
\def\fraku{{\frak u}}

\def\Fp{{{\Bbb F}_p}}
\def\GG{{\Bbb G}}
\def\NN{{\Bbb N}}
\def\QQ{{\Bbb Q}}
\def\ZZ{{\Bbb Z}}
\def\Zp{{ {\Bbb Z}_{(p)} }}

\def\calC{{\cal C}}
\def\calE{{\cal E}}
\def\calF{{\cal F}}
\def\calI{{\cal I}}
\def\calJ{{\cal J}}
\def\calL{{\cal L}}
\def\pcalL{\mbox{$p$-}{\cal L}}
\def\calM{{\cal M}}
\def\calP{{\cal P}}
\def\calS{{\cal S}}
\def\calU{{\cal U}}

\def\calAs{{\cal A}\hspace{-1pt}{\it s}\hspace{.3pt}}
\def\calLie{{\cal L}\hspace{-1pt}{\it ie}}
\def\pcalLie{\mbox{$p$-}{\cal L}\hspace{-1pt}{\it ie}}
\def\grcalLie{\mbox{gr-}{\cal L}\hspace{-1pt}{\it ie}}
\def\pgrcalLie{\mbox{$p$-gr-}{\cal L}\hspace{-1pt}{\it ie}}

\def\bV{{\bf V}}
\def\bW{{\bf W}}

\def\barCp{{\overline{C}_p}}
\def\barchi{{\overline{\chi}}}
\def\barFp{{\overline{\Bbb F}_p}}

\def\hatI{{\widehat{I}}}
\def\hatM{{\widehat{M}}}
\def\hatS{{\widehat{S}}}

\begin{document}

\title{Bernstein-Gelfand-Gelfand complexes and \\ 
cohomology of nilpotent groups over ${\bf Z}_{(p)}$ \\ 
for representations with $p$-small weights}
\author{P. Polo and J. Tilouine}
\date{\today}
\maketitle

{\hfuzz=3pt

\section*{Introduction}

Let $G$ be a connected reductive linear algebraic group defined and split over $\ZZ$, 
let $T$ be a maximal torus, $W$ the Weyl group, 
$R$ the root system, $R^\vee$ the set of coroots, $R^+$ a set of positive roots, 
and $\rho$ the half-sum of the elements of $R^+$. 
Let $X = X(T)$ be the character group of $T$ and let $X^+$ be the set of those 
$\lambda\in X$ such that $\langle\lambda, \alpha^\vee\rangle \geq 0$ for all $\alpha\in R^+$.

For any $\lambda\in X^+$, let $V_\ZZ(\lambda)$ be the Weyl module for $G$ over $\ZZ$ 
with highest weight $\lambda$ (see \ref{Weyl}) and, for any commutative ring $A$, let 
$V_A(\lambda) = V_\ZZ(\lambda)\otimes_\ZZ A$.   

Let $p$ be a prime integer and let  
$$
\barCp := \{\nu\in X \mid 0\leq \langle \nu+\rho,\beta^\vee\rangle \leq p,\quad 
\forall \beta\in R^+\},
$$
the closure of the fundamental $p$-alcove. 

The aim of this paper is to prove that several results about 
$V_\QQ(\lambda)$, due to Kostant \cite{Ko1}, 
Bernstein-Gelfand-Gelfand \cite{BGG}, Lepowsky \cite{Le}, Rocha \cite{Rocha}, and Pickel \cite{Pic},  
hold true over $\Zp$ when $\lambda\in X^+\cap \barCp$ : 
this is the precise meaning of the notion of $p$-smallness mentioned in the title.

\smallskip 
In more details, let $B$ be the Borel subgroup corresponding to $R^+$, let 
$P$ be a standard parabolic subgroup containing $B$, let $P^-$ be the opposed 
parabolic subgroup containing $T$, let $U_P^-$ be its unipotent radical, and let 
$L = P \cap P^-$, a Levi subgroup. Let $R_L$ the root system of $L$, let 
$R_L^+ = R_L\cap R^+$, and 
$$
X_L^+ := \{\xi\in X\mid \langle\xi,\alpha\rangle \geq 0 \quad\forall \alpha\in R_L^+\}.
$$
For any $\xi\in X_L^+$ and any commutative ring $A$, let $V_A^L(\xi)$ be the Weyl module 
for $L$ over $A$ with highest weight $\xi$.

Let $\frakg, \frakp, \fraku_P^-$ be the Lie algebras over $\ZZ$ of $G,P,U_P^-$, respectively, 
and let $U(\frakg)$ and $U(\frakp)$ be the enveloping algebras of $\frakg$ and $\frakp$.  
For $\xi\in X_L^+$, consider the generalized Verma module 
$$
M_\frakp^\ZZ(\xi) := U(\frakg)\otimes_{U(\frakp)} V_\ZZ^L(\xi).  
$$
For any commutative ring $A$, let $M_\frakp^A(\xi) = M_\frakp^\ZZ(\xi)\otimes_\ZZ A$. 
  
Let $N = \vert R^+\vert $ and, for $i= 0,1,\dots,N$, let 
$W(i) := \{w\in W \mid \ell(w) = i \}$, 
where $\ell$ denotes the length function on $W$ relative to $B$. 
Further, let $W^L = \{w\in W \mid wX^+ \subseteq X_L^+\}$ and 
$W^L(i) := W^L \cap W(i)$. 

\smallskip 
After several recollections in Section 1, we prove in Section the following Theorem 
(see \ref{weakBGG}).  

\begin{thintro}\label{introweakBGG} Let $\lambda\in X^+\cap \barCp$. 
There exists an exact sequence of $U(\frakg)$-modules,  
$$
0\to D_N(\lambda) \to \cdots \to D_0(\lambda) \to V_\Zp(\lambda) \to 0,
$$ 
where each $D_i(\lambda)$ admits a finite filtration of $U(\frakg)$-submodules with associated 
graded 
$$
\gr D_i(\lambda) \cong \bigoplus_{w\in W^L(i)} M_\frakp^\Zp(w(\lambda+\rho)-\rho).
$$
\end{thintro} 

That is, following the terminology introduced in \cite{Rocha}, 
$V_\Zp(\lambda)$ admits a weak generalized 
Bernstein-Gelfand-Gelfand resolution. From this, one obtains immediately 
the following (see \ref{kostant} and \ref{pf-kostant}). 

\begin{thintro}\label{introKostant}{\rm (Kostant's theorem over $\Zp$)} 
Let $\lambda\in X^+\cap \barCp$. Then, 
for each $i$, there is an isomorphism of $L$-modules : 
$$
H_i(\fraku_P^-, V_\Zp(\lambda)) \cong \bigoplus_{w\in W^L(i)} V_\Zp^L(w(\lambda+\rho)-\rho).
$$
\end{thintro}

 Let $\Gamma := U_P^-(\ZZ)$ be the 
group of $\ZZ$-points of $U_P^-$, it is a finitely generated, torsion free, 
nilpotent group. By a result of Pickel \cite{Pic}, there is a natural  
isomorphism $H_*(\fraku_P^-, V_\QQ(\lambda)) \cong H_*(\Gamma, V_\QQ(\lambda))$. 
In Section 3, we prove a slightly weaker version of this result 
over $\Zp$ when $\lambda$ is $p$-small (see \ref{mainthree}).   

\begin{thintro}\label{introPickel} Let $\lambda\in X^+\cap \barCp$. Then, for each $n\geq 0$,  
$H_n(U_P^-(\ZZ), V_\Zp(\lambda))$ has a natural $L(\ZZ)$-module filtration 
such that 
$$
gr H_n(U_P^-(\ZZ), V_\Zp(\lambda)) \cong \bigoplus_{w\in W^L(n)} V_\Zp^L(w(\lambda+\rho)-\rho).
$$
\end{thintro} 

The proof of this result has two parts. In the first, we develop certain general results 
about finitely generated, torsion free, nilpotent group $\Gamma$. In particular, using a 
beautiful theorem of Hartley \cite{Hart}, which perhaps did not receive as much attention 
as it deserved, we obtain in an algebraic manner a spectral sequence 
relating the homology of a certain graded, torsion-free, Lie ring 
$\grisol\Gamma$ associated with $\Gamma$ 
to the homology of $\Gamma$ itself, the coefficients being a $\Gamma$-module 
with a ``nilpotent" filtration and its associated graded 
(see Theorem \ref{twospecseq}). 
This gives a purely algebraic, homological version (with coefficients) of a 
cohomological spectral sequence obtained, using methods of algebraic topology, 
by Cenkl and Porter \cite{CP}. 
In fact, our methods also have a cohomological counterpart. 
This will be developped in a subsequent paper \cite{Po}. 

In the second part of the proof, we first show that in our case where $\Gamma = U_P^-(\ZZ)$, 
one has $\grisol \Gamma \cong \fraku_P^-$, and then deduce from the truth of Kostant's theorem 
over $\Zp$ that the spectral sequence mentioned above degenerates at $E_1$.

Next, in Section 4, we obtain a result \`a la Bernstein-Gelfand-Gelfand concerning 
now the distribution algebras $\dist(G)$ and $\dist(P)$. 
Namely, there exists a standard {\it complex\,} (not a resolution!)
$$
\calS_\bullet(G,P,\lambda) = 
\dist(G)\otimes_{\dist(P)} \left(\Lambda^\bullet(\frakg/\frakp) \otimes V_\ZZ(\lambda)\right).  
$$
For $\xi\in X_L^+$, consider the generalized Verma module (for $\dist(G)$ and $\dist(P$)) 
$$
\calM_P^\ZZ(\xi) := \dist(G)\otimes_{\dist(P)} V_\ZZ^L(\xi),  
$$
and, for any commutative ring $A$, set $\calM_P^A(\xi) = \calM_P^\ZZ(\xi)\otimes_\ZZ A$. 

 Under the assumption that the unipotent radical of $P$ 
is abelian, we obtain, by using an idea borrowed from 
\cite[\S\,VI.5]{CF} plus arguments from 
Section 2, the following result (see \ref{th-BGGcomplex}).

\begin{thintro}\label{introBGGdist} Suppose that $\fraku_P^-$ is abelian. 
Let $\lambda\in X^+\cap\barCp$. Then 
the standard complex $\calS^\Zp_\bullet(G,P,\lambda)$ contains as a direct summand 
a subcomplex $\calC^\Zp_\bullet(G,P,\lambda)$ such that, for $i\geq 0$, 
$$
\calC_i^\Zp(G,P,\lambda) 
\cong 
\bigoplus_{w\in W^L(i)} \calM_P^\Zp(w\cdot\lambda).
$$
\end{thintro}

Presumably, the hypothesis that $U_P^-$ be abelian can be removed, but the proof 
would then require considerably more work. Since the  abelian case is 
sufficient for the applications in the companion 
paper by A. Mokrane and J. Tilouine \cite{MT}, we content ourselves with this 
result. We hope to come back to the general case later.

\medskip 
To conclude this introduction, let us mention that the results of this text 
are used in \cite{MT} in the case where $G$ is the group of symplectic similitudes.   
When $P$ is the Siegel parabolic, Theorem \ref{introBGGdist} occurs in \cite[\S\,5.4]{MT} as 
an important step to establish a modulo $p$
analogue of the Bernstein-Gelfand-Gelfand complex of 
\cite[Chap.VI, Th.\,5.5]{CF}, while Theorem \ref{introPickel} (in its cohomological form) 
is used in \cite[\S\,8.3]{MT} to study mod. $p$ versions of Pink's theorem on higher 
direct images of automorphic bundles.

The notations of \cite{MT}
follow those of \cite{CF} and are therefore different from the ones used in the present paper, which
are standard in the theory of  reductive groups. 
A dictionary is provided in the final section  of this text.

\section{Notation and preliminaries}

\subsection{}  
Let $G$ be a connected reductive linear algebraic group, defined and split over $\ZZ$. 
Let $T$ be a maximal torus, $W$ the Weyl group, $R$ the root system and $R^\vee$ the set of coroots. 
Fix a set $\Delta$ of simple roots, let $R^+$ and $R^-$ be the corresponding sets of positive and 
negative roots, and let $B$, $B^-$ denote the associated Borel subgroups and $U$, $U^-$ 
their unipotent radicals. (For all this, see, for example, \cite{Dem}). 

Let $X = X(T)$ be the character group of $T$; elements of $X$ will be called weights, in accordance 
with the terminology in Lie theory. 
Let $\leq$ denote the partial order on $X$ defined by the positive cone $\NN R^+$, that is, 
$\mu\leq \lambda$ if and only if $\lambda-\mu \in \NN R^+$. 
Let $Q(R)\subset X$ be the root lattice and let $\rho$ 
be the half-sum of the positive roots; it belongs to $X\otimes \ZZ[1/2]$. Define, as usual, the 
dot action of $W$ on $X$ by 
$$
w\cdot \lambda = w(\lambda+\rho) - \rho,
$$
for $\lambda\in X, w\in W$. It is well-known and easy to see that $w\rho-\rho\in Q(R)$ for 
all $w\in W$, and hence $w\cdot \lambda$ does indeed belong to $X$. 

Let $X^+$ be the set of dominant weights: 
$$
X^+ := \{\lambda\in X\mid \forall\alpha\in R^+, \quad \langle\lambda,\alpha^\vee\rangle \geq 0\},
$$
where $\alpha^\vee$ denotes the coroot associated with $\alpha$.

\subsection{Enveloping and distribution algebras.}\label{PvsDist}
Let $\frakg = \Lie(G)$ (resp. $\frakt = \Lie(T)$) be the Lie algebra of $G$ (resp. $T$); 
they are finite free $\ZZ$-modules. 
Let $U(\frakg)$ denote the enveloping algebra of $\frakg$ over $\ZZ$, and let 
$\dist(G)$ denote the algebra of distributions of $G$ 
(see \cite[Chap.\,I.7]{Ja}). 
If $G$ is semi-simple 
and simply-connected, $\dist(G)$ coincides with the Kostant $\ZZ$-form of $U(\frakg)$ 
(\cite{Ko2}), 
see \cite[\S\,II.1.12]{Ja} or \cite[VIII, \S\S\,12.6--8]{Bki}. 
We shall denote it by $\calU_\ZZ(\frakg)$ or simply $\calU(\frakg)$; sometimes it will also be 
convenient to denote it by $\calU_\ZZ(G)$.  

Similarly, if $H$ is a closed subgroup of $G$ defined over $\ZZ$, we shall denote 
$\dist(H)$ also by $\calU_\ZZ(H)$.

By an $H$-module we shall mean a rational $H$-module, that is, a $\ZZ[H]$-comodule. 
More generally, for any commutative ring $A$, an $H_A$-module means an 
$A$-module with a structure of $A[H]$-comodule. 
If $V$ is an $H$-module, then, as is well-known, $V$ is also 
an $\calU_\ZZ(H)$-module and a fortiori an $U(\Lie(H))$-module. 

If $M$ is a $T$-module, it is the direct sum of its weight spaces $M_\lambda$, 
for $\lambda\in X$, see, for example, \cite[\S\,I.2.11]{Ja}. 

For future use, let us record here the following 

\begin{prop} Let $P$ be a standard parabolic subgroup of $G$, 
let $V$ be a finite dimensional $P_\QQ$-module and let $M$ be a 
$\ZZ$-lattice in $V$. Then $M$ is a $P$-submodule if and only is it is 
an $\calU_\ZZ(P)$-submodule. 
\end{prop}

\pf Without loss of generality we may assume that $P$ contains $B$. 
Let $P^-$ be the opposed standard parabolic subgroup 
and let $U_P^-$ be its unipotent radical. 
By the Bruhat decomposition, the multiplication map 
induces an isomorphism of $U_P^- \times B$ onto an open subset of $P$, 
see, for example, \cite[\S\,II.1.10]{Ja}. 
This implies that the arguments in \cite[II.8.1]{Ja} 
are valid for $P$, and the proposition then follows from \cite[I.10.13]{Ja}. 

\subsection{Weyl modules.}\label{Weyl} 

For $\lambda\in X^+$, let $V_\QQ(\lambda)$ denote the irreducible $G_\QQ$-module with highest 
weight $\lambda$, and let $V_\ZZ(\lambda)$ be the corresponding Weyl module over $\ZZ$; that is, 
$$
V_\ZZ(\lambda) := \calU_\ZZ(G) v_\lambda 
$$
is the $\calU_\ZZ(G)$-submodule generated by a fixed vector $v_\lambda \not= 0$ of weight $\lambda$. 
It is a $G$-module by Proposition \ref{PvsDist} above. 
Of course, up to isomorphism, $V_\ZZ(\lambda)$ does not depend on the choice of $v_\lambda$. 
For future use, let us also record the following (obvious) lemma. 

\begin{lem} Let $M$ be a $\ZZ$-free $G$-module and $v\in M$ an element fixed by $U$ and of weight 
$\lambda$. Then the submodule $\calU_\ZZ(G) v$ is isomorphic to $V_\ZZ(\lambda)$. 
\end{lem}

\pf The $\calU_\QQ(G)$-submodule of $M\otimes \QQ$ generated by $v$ is isomorphic to 
$V_\QQ(\lambda)$.

\subsection{Contravariant duals.}\label{def-tau} 
Let us fix an anti-involution $\tau$ of $G$ which is the identity on $T$ and 
exchanges $B$ and $B^-$ (see \cite[II.1.16]{Ja}). Then $\tau$ induces 
anti-involutions on $\calU_\ZZ(G)$, on $\frakg$ and on $U_\ZZ(\frakg)$, 
which we denote by the same letter 
$\tau$. 

For any ring $A$ and $G_A$-module $V$, let us denote by $V^\tau$ the $A$-dual 
$\Hom_A(V,A)$, regarded as a $G_A$-module via $\tau$. It may be called the 
``contravariant dual" of $V$, as for $V = V_\ZZ(\lambda)$ this is closely related 
to the so-called ``contravariant form" on $V_\ZZ(\lambda)$; see \cite[II.8.17]{Ja} 
and the discussion in the next subsection \ref{lattices}. 

Note that if $V$ is a free $A$-module, the weights of $T$ in $V$ and $V^\tau$ are the 
same. In particular, the irreducible $G_\QQ$-modules $V_\QQ(\lambda)$ and 
$V_\QQ(\lambda)^\tau$ are isomorphic.

\subsection{Admissible lattices.}\label{lattices}
For use in the companion article by Mokrane and Tilouine \cite{MT}   
 and also in the next subsection, let us discuss some properties of 
admissible lattices. Of course, this is fairly well-known to representation theoretists, 
but we spell out the details for the convenience of readers with a different background. 

As noted above, we may identify $V_\QQ(\lambda) = V_\QQ(\lambda)^\tau$. Under this 
identification, $V_\QQ(\lambda)$ becomes equipped with a 
non-degenerate, $G$-invariant bilinear form 
$\langle\phantom{v},\phantom{v}\rangle$ such that 
$$
\langle gv, v'\rangle = \langle v, \tau(g)v'\rangle \quad\mbox{and}\quad 
\langle Xv, v'\rangle = \langle v, \tau(X)v'\rangle, \leqno (*)
$$
for $v,v'\in V_\QQ(\lambda)$, $g\in G$, $X\in \calU_\ZZ(G)$. (This is the contravariant form 
mentioned in the previous subsection). 

Let us fix, once for all, a non-zero vector $v_\lambda\in V_\QQ(\lambda)_\lambda$. 
The identification $V_\QQ(\lambda) = V_\QQ(\lambda)^\tau$ may be chosen so that 
$\langle v_\lambda, v_\lambda\rangle = 1$.  
  
Recall that a $\ZZ$-lattice $\calL \subset V_\QQ(\lambda)$ is called an admissible lattice 
if it is stable under $\calU_\ZZ(G)$. By Proposition \ref{PvsDist}, this implies that 
$\calL$ is a $G$-module and is therefore the direct sum of its $T$-weight spaces. 

Let $\calE(\lambda)$ denote the set of admissible lattices $\calL\subset V_\QQ(\lambda)$ 
such that $\calL\cap V_\QQ(\lambda)_\lambda = \ZZ v_\lambda$. Clearly, 
$V_\ZZ(\lambda) := \calU_\ZZ(G) v_\lambda$ is the unique minimal element of $\calE(\lambda)$. 

For any $\calL\in \calE(\lambda)$, the dual $G$-module $\calL^\tau$ identifies with 
$$
\{x\in V_\QQ(\lambda) \mid \langle x, \calL\rangle \subseteq \ZZ\}.
$$
It follows from $(*)$ that $\calL^\tau$ is an admissible lattice, 
and since $\langle v_\lambda, v_\lambda\rangle = 1$ it belongs to $\calE(\lambda)$. 
Therefore, $\calL^\tau \supseteq V_\ZZ(\lambda)$ and hence 
$\calL \subseteq V_\ZZ(\lambda)^\tau$. Let us record this as the next 

\begin{lem} The set of admissible lattices $\calL\subset V_\QQ(\lambda)$ 
such that $\calL\cap V_\QQ(\lambda)_\lambda = \ZZ v_\lambda$ contains a 
unique minimal element, $V_\ZZ(\lambda)$, and a unique maximal element, 
$V_\ZZ(\lambda)^\tau$.
\end{lem} 

The above minimal and maximal lattices are denoted by $V(\lambda)_{min}$ and 
$V(\lambda)_{max}$ in \cite{MT} and in Section \ref{dict} below. 

\subsection{Weyl modules and induced modules.}\label{Weyl-ind}

Let us recall the definition of the induction functor $\Ind_{B^-}^G : 
\{B^-\mbox{-modules}\} \to \{G\mbox{-modules}\}$. For any $B^-$-module $M$, 
$$
\Ind_{B^-}^G(M) := (\ZZ[G]\otimes M)^{B^-},
$$
where $\ZZ[G]$ is regarded as a $G\times B^-$-module via 
$\left((g,b)\phi\right)(g') = \phi(g^{-1}g'b)$, for 
$g,g'\in G$, $b\in B^-$ and where the invariants are taken with respect to the diagonal 
action of $B^-$; it is a left exact functor, see \cite[\S\,I.3.3]{Ja}. 
As in \cite[\S\,II.2.1]{Ja}, we shall denote simply by $H^i(\phantom{-})$ the 
right derived functors $R^i\, \Ind_{B^-}^G(\phantom{-})$. 

Let $\lambda\in X$; it may be regarded in a natural manner as a character of either 
$B^-$ or $B$. Moreover, since $\tau$ is the identity on $T$, one has 
$\lambda(\tau(b)) = \lambda(b)$ for any $b\in B^-$. 

For any ring $A$, let us denote by $A_\lambda$ the free $A$-module 
of rank one on which $B^-$ acts via the character $\lambda$. Then, 
$$
H^0(A_\lambda) \cong \{\phi\in A[G] \mid \phi(gb) = \lambda(b^{-1})\phi(g), \;\; 
\forall\, g\in G, b\in B^-\}.
$$

\begin{prop} Let $\lambda\in X^+$. 

{\rm a)} $H^0(\ZZ_\lambda) \cong V_\ZZ(\lambda)^\tau$. 

{\rm b)} If $k$ is a field, 
$H^0(k_\lambda) \cong H^0(\ZZ_\lambda)\otimes k \cong V_k(\lambda)^\tau$. 
Thus, in particular, $V_k(\lambda)$ is irreducible if and only if 
$H^0(k_\lambda)$ is so. 
\end{prop}

\pf First, by flat base change (\cite[I.3.5]{Ja}), one has $H^0(\ZZ_\lambda)\otimes \QQ 
\cong H^0(\QQ_\lambda)$. Moreover, $H^0(\QQ_\lambda) \cong V_\QQ(\lambda)$, by the 
theorem of Borel-Weil-Bott (see, for example, \cite[II.5.6]{Ja}).

Further, since $\ZZ[G]$ is a free $\ZZ$-module (being a subring of 
$\ZZ[U]\otimes \ZZ[B^-]$), so is 
$H^0(\ZZ_\lambda)$. Therefore, $H^0(\ZZ_\lambda)$ may be identified with a 
$G$-submodule of $V_\QQ(\lambda)$, and the identification may be chosen so that 
$H^0(\ZZ_\lambda) \cap V_\QQ(\lambda)_\lambda = \ZZ v_\lambda$, {\it i.e.}, 
so that $H^0(\ZZ_\lambda)$ belongs to $\calE(\lambda)$. 

Now, there is a natural $G$-module map $\phi : V_\ZZ(\lambda)^\tau \to H^0(\ZZ_\lambda)$ 
defined by 
$$
x\mapsto \left( g\mapsto \langle x, \tau(g^{-1})v_\lambda \rangle \right).
$$
Moreover, since $V_\ZZ(\lambda)$ is generated by $v_\lambda$ as a $G$-module, 
$\phi$ is injective. Since $V_\ZZ(\lambda)^\tau$ is the largest element of $\calE(\lambda)$, 
this implies that $\phi$ induces an isomorphism $V_\ZZ(\lambda)^\tau \cong H^0(\ZZ_\lambda)$. 
This proves assertion a). 

\smallskip Let us prove assertion b). For each $i\geq 0$ there is an exact sequence 
$$
0\to H^i(\ZZ_\lambda)\otimes k\to H^i(k_\lambda)\to \Tor^\ZZ(H^{i+1}(\ZZ_\lambda),k)\to 0,
$$
see \cite[I.4.18]{Ja}. Next, by Kempf's vanishing theorem (\cite[II.4.6]{Ja}), one has 
$H^i(\ZZ_\lambda) = 0$ for $i\geq 1$. The first isomorphism of assertion b) follows. 
Finally, the second is a consequence of assertion a) and the natural isomorphisms 
$$
\Hom_\ZZ(V_\ZZ(\lambda),\ZZ)\otimes k \cong \Hom_\ZZ(V_\ZZ(\lambda),k) \cong 
\Hom_k(V_k(\lambda),k).
$$
This completes the proof of the proposition.

\subsection{Parabolic subgroups and unipotent radicals.} 
Now, let $P$ be a standard parabolic subgroup of $G$ containing $B$, let $L$ 
be the Levi subgroup of $P$ containing $T$, and let $P^-$ be the standard parabolic 
subgroup opposed 
to $P$, that is, $P^-$ is the unique parabolic subgroup containing $B^-$ such that 
$P^- \cap P = L$. 

Let $U_P^-$ (resp.\,$U_P$) denote the unipotent radical of $P^-$ (resp.\,$P$) 
and let $\fraku_P^- = \Lie(U_P^-)$, $\fraku_P = \Lie(U_P)$ and 
$\frakp = \Lie(P)$. Then $\fraku_P^-$, $\fraku_P$ and $\frakp$ are free $\ZZ$-modules and 
$\frakg = \frakp \oplus \fraku_P^- $.
Thus, in particular, $\frakg/\frakp$ is a free $\ZZ$-module. 

Further, if $V$ is a $P$-module then, by standard arguments, the homology groups 
$$
H_i(\fraku_P^-, V) := \Tor_i^{U(\fraku_P^-)}(\ZZ,V)
$$
carry a natural structure of $L$-modules. For example, they can be computed as the homology of 
the standard Chevalley-Eilenberg complex $\Lambda^\bullet(\fraku_P^-)\otimes V$, which 
carries a natural action of $L$. 

For any commutative ring $A$, we set $V_A(\lambda) := V_\ZZ(\lambda)\otimes A$ 
and $\frakg_A := \frakg\otimes A$. The enveloping algebra of $\frakg_A$ identifies with 
$U_\ZZ(\frakg)\otimes A$ and is denoted by $U_A(\frakg)$. One defines similarly $U_A(\fraku_P^-)$ 
and $\calU_A(\frakg)$, etc...  

Since $U_\ZZ(\fraku_P^-)$ is a free $\ZZ$-module, one has, for every $i\geq 0$, 
$$
\Tor_i^{U_A(\fraku_P^-)}(A, V_A(\lambda)) \cong \Tor_i^{U_\ZZ(\fraku_P^-)}(\ZZ, V_A(\lambda)).
$$
We shall denote these groups simply by $H_i(\fraku_P^-,V_A(\lambda))$; as noted above they are 
$L_A$-modules.

Our goal in Section 2 is to show that celebrated results of Kostant (\cite[Cor.\,8.1]{Ko1}) 
and Bernstein-Gelfand-Gelfand (\cite[Th.\,9.9]{BGG}), which describe respectively, 
for any $\lambda\in X^+$,   
the $L$-module structure of $H_\bullet(\fraku_P^-, V_\QQ(\lambda))$ 
and a minimal $U_\QQ(\fraku_P^-)$-resolution of $V_\QQ(\lambda)$, 
hold true 
when $\QQ$ is replaced by $\ZZ_{(p)}$, for any prime integer $p$ such that 
$$
p \geq \langle\lambda+\rho, \alpha^\vee\rangle, \quad \forall\,\alpha\in R^+ .
$$

\subsection{Weyl modules for a Levi subgroup.}\label{defWL}   
We need to introduce more notation. Let $W_L$ and $R_L$ denote the Weyl group and root system 
of $L$, and let $R_L^+ := R_L \cap R^+$. Let $X_L^+$ denote the set of 
$L$-dominant weights  
$$
X_L^+ := \{\lambda\in X\mid \forall\alpha\in R_L^+, \quad \langle\lambda,\alpha^\vee\rangle 
\geq 0\}.
$$ 

Let $W^L := \{w\in W\mid wX^+ \subseteq X_L^+\}$. It is well-known, and easy to check, 
that $W^L$ is also equal to $\{w\in W \mid w^{-1}R_L^+ \subseteq R^+\}$. 

Let $\ell$ and $\leq$ denote the length function and Bruhat-Chevalley order on $W$ 
associated with the set $\Delta$ of simple roots. 
Then, for $i\geq 0$, set 
$$
W(i) := \{w\in W\mid \ell(w) = i\} \quad \mbox{and}\quad W^L(i) := W^L\cap W(i).
$$

For any $\xi\in X_L^+$, let $V_\QQ^L(\xi)$ denote the irreducible $L_\QQ$-module 
with highest weight $\xi$ and let $V_\ZZ^L(\xi)$ be the corresponding Weyl module 
for $L$. Observe that $V_\QQ^L(\xi)$ (and then $V_\ZZ^L(\xi)$) identifies with the 
$L_\QQ$-submodule of $V_\QQ(\xi)$ (resp. $L$-submodule of $V_\ZZ(\xi)$) generated by 
$v_\xi$. 

More generally, one has the following 

\begin{lem} Let $M$ be a $P$-module which is $\ZZ$-free and let $v\in M$ be a non-zero element  
of weight $\xi$. Assume that $v$ is $U$-invariant $($this is the case, for instance, if $\xi$ is 
a maximal weight of $M)$. Then the $\calU_\ZZ(P)$-submodule of $M$ generated by $v$ is isomorphic to 
$V_\ZZ^L(\xi)$. 
\end{lem} 

\pf Recall that 
$\calU_\ZZ(P) \cong \calU_\ZZ(L)\otimes \calU_\ZZ(U_P)$ (see \cite[\S\,II.1.12]{Ja}).
Since $v$ is fixed by $U$, it is annihilated by the augmentation ideal of 
$\calU_\ZZ(U_P)$. Therefore, $\calU_\ZZ(P)v = \calU_\ZZ(L) v$ and, since $M$ 
is $\ZZ$-free, the result follows from Lemma \ref{Weyl}.

\subsection{The fundamental $p$-alcove.}\label{alcove}\label{irred}  
In this subsection and the next one, let $p$ be a prime integer. 
The notion of $p$-smallness mentioned in the title of this article 
is defined as follows. We shall say that $\lambda\in X$ is $p$-small if it satisfies the 
condition:
$$
 \langle\lambda+\rho,\alpha^\vee\rangle \leq p,\quad \forall\,\alpha\in R .\leqno (\dagger)
$$

\smallskip 
An equivalent definition of $p$-smallness is as follows. 
Let $W_p$ denote the affine Weyl group with respect to 
$p$. Recall that $W_p$ is the subgroup of automorphisms of $X(T)\otimes \Bbb R$ generated by the 
reflections $s_{\beta, np}$, for $\beta\in R^+$, $n\in \ZZ$, where, 
for $\lambda\in X(T)\otimes \Bbb R$,  
$$
s_{\beta, np}(\lambda) = \lambda - (\langle\lambda,\beta^\vee\rangle - n)\beta, 
$$
and that $W_p$ is the semi-direct product of 
$W$ and the group $pQ(R)$ acting by translations.  We consider the dot action of $W_p$ on 
$X(T)\otimes \Bbb R$, defined by $w\cdot\lambda = w(\lambda+\rho)-\rho$. 

The fundamental $p$-alcove $C_p$ is defined by 
$$
C_p := \{\lambda\in X(T)\otimes {\Bbb R} \mid 0<\langle\lambda+\rho,\beta^\vee\rangle < p, 
\;\; \forall\,\beta\in R^+\}.
$$
Its closure  
$$
\barCp := \{\lambda\in X(T)\otimes {\Bbb R} \mid 
0\leq \langle\lambda+\rho,\beta^\vee\rangle \leq p, 
\;\; \forall\,\beta\in R^+\}
$$ 
is a fundamental domain for the dot action of $W_p$ on $X(T)\otimes \Bbb R$ 
(for all this, see for example \cite[\S\,II.6.1]{Ja}). 

Then, for $\lambda\in X^+$, the condition of $p$-smallness is equivalent to the requirement 
that $\lambda$ belongs to $\barCp$. 

Let $\rho_L$ be the half-sum of the elements of $R_L^+$. Note that 
$\langle \rho_L, \alpha^\vee\rangle = 1$ for any $\alpha\in \Delta\cap R_L$ and hence 
$\rho-\rho_L$ vanishes on $R_L$. Therefore, if a weight $\xi\in X_L^+$ is $p$-small, 
it is a fortiori $p$-small for $L$. 

\medskip The fact that $V_\Fp(\lambda)$ is irreducible when $\lambda$ is $p$-small is of course 
very well-known to representation-theoretists; for the convenience of readers with a different 
background, we record this here as the next 

\begin{lem} Let $\lambda\in X^+$ and $\xi\in X_L^+$. 
If $\lambda$ (resp. $\xi$) is $p$-small, $V_\Fp(\lambda)$ (resp. $V_\Fp^L(\xi)^\tau$) 
is irreducible 
and self-dual for the contravariant duality.  
\end{lem}

\pf 
The first assertion is a consequence of \cite[II.8.3]{Ja}, 
combined with Proposition \ref{Weyl-ind}. 
Further, since irreducible $G_\Fp$-modules are determined by their highest weight, 
the second assertion follows from the first. 

\begin{cor} If $\lambda\in X^+\cap\barCp$ then, for any $\Lambda\in\calE(\lambda)$, one has 
$$
V_\Zp(\lambda) = \Lambda\otimes \Zp = V_\Zp(\lambda)^\tau.
$$
\end{cor}

\pf By the previous lemma, one has $V_\Fp(\lambda) = V_\Fp(\lambda)^\tau$.   
The result then follows by Nakayama's lemma.

\subsection{A vanishing result.}\label{vanish}
Let us record the following

\begin{lem} For all $\lambda,\mu \in X^+$, one has $\Ext^1_G(V_\Fp(\lambda),V_\Fp(\mu)^\tau)= 0$ 
and also 
$$
\Ext^1_G(V_\ZZ(\lambda),V_\ZZ(\mu)^\tau)= 0 = \Ext^1_G(V_\Zp(\lambda),V_\Zp(\mu)^\tau).   
$$  
\end{lem}

\pf Since $V_\Fp(\mu)^\tau \cong H^0(\mu)$, by Proposition \ref{Weyl-ind},  
the assertion over $\Fp$ is a consequence of \cite[Prop.\,II.4.13]{Ja}. 
The assertions over $\ZZ$ or $\Zp$ then follow from a 
theorem of universal coefficients \cite[Prop.\,I.4.18]{Ja}.

\begin{cor} Suppose that $\lambda,\mu \in X^+\cap \barCp$. Then 
$$
\Ext^1_G(V_\Fp(\lambda),V_\Fp(\mu)) = 0 = \Ext^1_G(V_\Zp(\lambda),V_\Zp(\mu)). 
$$ 
\end{cor}

\pf By the results in \ref{irred}, $V_\Fp(\mu)$ and $V_\Zp(\mu)$ are self-dual. 
Thus, the corollary follows from the previous lemma.

\subsection{}\label{donkin} 
We shall need later the following 

\begin{lem} Let $M$ be a $P$-module, finite free over $\ZZ_{(p)}$. Assume that each weight 
$\nu$ of $M$ satisfies $\langle\nu+\rho, \alpha^\vee\rangle \leq p$, for any $\alpha\in R_L$. 
Then there exists a sequence of $P$-submodules 
$0 = M_0\subset \cdots \subset M_r = M$ such that 
$$
M_i/M_{i-1} \cong V_{\Zp}^L(\xi_i), \mbox{ for some } \xi_i\in X_L^+
$$
and $\xi_j\leq \xi_i$ if $j\geq i$. Further, the set $\{\xi_1,\dots,\xi_r\}$ is uniquely 
determined by $M$; in fact the $V_{\QQ}^L(\xi_i)$ are the irreducible composition factors 
of the $L_\QQ$-module $M_\QQ$. 
\end{lem} 

\pf (following \cite[Lemma 11.5.3]{Do}) We proceed by induction on $d = \dim_\QQ M_\QQ$. 
There is nothing to prove if 
$M = 0$. If $M\not= 0$, let $\xi_1$ be a maximal weight of $M$, let $v\in M$ be a primitive  
element of weight $\xi_1$ and denote by $N$ the $\calU_\Zp(P)$-submodule generated by $v$.  
Then $N \cong V_{\Zp}^L(\xi_1)$. By assumption, $\xi_1\in \barCp$ and hence 
$N_\Fp := N\otimes \Fp$ is irreducible. 

On the other hand, since $M$ is free over $\Zp$, one obtains an exact sequence of $P$-modules 
$$
0\to \Tor_1^\Zp(M/N,\Fp) \to N_\Fp\stackrel{\phi}{\to} M_\Fp ,
$$
and $\phi(v)\not= 0$, as $v$ is a primitive element. Since $N_\Fp$ is 
irreducible, $\phi$ is injective. Thus, $\Tor_1^\Zp(M/N,\Fp) = 0$ and it follows that $M/N$ 
is free over $\Zp$. The first assertion of the lemma then follows by the induction hypothesis.
Finally, the second assertion is clear.

\section{Kostant's theorem over ${\bf Z}_{(p)}$}

\subsection{}\label{cox}\label{kostant} 
 In this section, let us fix $\lambda\in X^+$ and let $p$ be a prime integer 
such that 
$$
p\geq \langle \lambda+\rho,\alpha^\vee\rangle, \quad \forall\,\alpha\in R^+ . \leqno (\dag) 
$$

{\it Remark}. It is customary, in representation theory, to introduce the so-called 
Coxeter number of $G$, defined by 
$$
h := 1 + \Max\{\langle\rho,\alpha^\vee\rangle,\; \alpha\in R^+\}. 
$$
Therefore, since $\lambda$ is dominant, our assumption $(\dag)$ above implies that 
$p\geq h-1$, and reduces to this inequality in the case where $\lambda = 0$. 
Thus, we allow the cases where $p$ equals $h-1$ or $h$. 

\medskip 
Our goal in this section is to prove the following 

\begin{th} Let $\lambda\in X^+$ and $p$ as above. Then, for each $i$, there is an 
isomorphism of $L$-modules 
$$
H_i\bigl(\fraku_P^-, V_{\ZZ_{(p)}}(\lambda)\bigr) \cong \bigoplus_{w\in W^L(i)} 
V_{\ZZ_{(p)}}^L(w\cdot\lambda).
$$
\end{th}

By standard arguments, it suffices to prove the theorem in the case where $G$ is semi-simple 
and simply-connected and the root system $R$ is irreducible. Similarly, the result for 
$SL_n$ is easily derived from the result for $GL_n$ (for technical reasons, the latter is 
easier to handle, see below). Therefore, we may (and shall) assume in the rest of this section that 
$R$ is irreducible and that $G$ is either $GL_n$ or semi-simple and simply-connected of type 
$\not =A$.

\subsection{Standard resolutions for $U({\frak g})$.}\label{standard} 
Recall first the standard Koszul resolution of the trivial 
module:
$$
\cdots \to U(\frakg)\otimes \Lambda^2(\frakg) \stackrel{d_2}{\lra} 
U(\frakg)\otimes \frakg \stackrel{d_1}{\lra} U(\frakg) \stackrel{\varepsilon}{\lra} \ZZ \to 0,
$$
where each differential $d_k$ is defined by the formula 
$$
\displaylines{ 
d_k(u\otimes x_1\wedge\cdots\wedge x_k) := \sum_{i=1}^k (-1)^{i-1}\, u x_i\otimes 
x_1\wedge\cdots\wedge \widehat{x_i}\wedge \cdots \wedge x_k 
\hfill\cr\hfill  
+ \sum_{1\leq i<j\leq k} (-1)^{i+j}\, u\otimes [x_i,x_j]\wedge x_1\wedge 
\cdots\wedge \widehat{x_i}\wedge \cdots \widehat{x_j}\wedge \cdots\wedge x_k. 
}
$$
Let $\pi_\frakp$ denote the natural projection 
$\Lambda^\bullet(\frakg)\to \Lambda^\bullet(\frakg/\frakp)$; it is a morphism 
of $P$-modules. Then, there is a surjective 
morphism of $U(\frakg)$-modules: 
$$
\begin{array}{cc}
\phi_\frakp :& U(\frakg)\otimes \Lambda^\bullet(\frakg) \to U(\frakg)\otimes_{U(\frakp)} 
\Lambda^\bullet(\frakg/ \frakp) \\
& u\otimes x \mapsto u\otimes_{U(\frakp)} \pi_\frakp(x).
\end{array}
$$
It is well-known, and easy to check, that each $d_k$ induces a map 
$d_k^\frakp$ such that $\phi_\frakp \circ d_k = d_k^\frakp \circ \phi_\frakp$.  
Thus, one obtains a complex of $U(\frakg)$-modules 
$$
\cdots \to U(\frakg)\otimes_{U(\frakp)} \Lambda^2(\frakg/ \frakp) 
 \stackrel{d_2^\frakp}{\to} U(\frakg)\otimes_{U(\frakp)} (\frakg/ \frakp)   
\stackrel{d_1^\frakp}{\to} U(\frakg)\otimes_{U(\frakp)} \ZZ 
\stackrel{\varepsilon}{\to} \ZZ\to 0,  
$$
which is still exact, for it is easily seen that the proof of \cite[Th.9.1]{BGG} 
is valid over $\ZZ$. This complex is called the standard resolution 
of the trivial module $\ZZ$ relative to $U(\frakg)$ and $U(\frakp)$. We shall denote it by 
$S_\bullet(\frakg,\frakp,\ZZ)$ or simply $S_\bullet(\frakg,\frakp)$. 

\smallskip 
Let $V$ be a $\ZZ$-free $U(\frakg)$-module. Then 
$S_\bullet(\frakg,\frakp)\otimes V$, with the diagonal 
action of $\frakg$, is an $U(\frakg)$-resolution of $V$ by modules which are free over 
$U(\fraku_P^-)$. 

Further, recall the ``tensor identity" \cite[Prop.\,1.7]{GL} : 
for any $U(\frakp)$-module $E$, there is 
a natural isomorphism of $U(\frakg)$-modules 
$$
(U(\frakg)\otimes_{U(\frakp)} E)\otimes V \cong 
U(\frakg)\otimes_{U(\frakp)} (E\otimes V_\vert{}_\frakp),
$$
where $V_\vert{}_\frakp$ denotes $V$ regarded as an $U(\frakp)$-module. Applying these 
isomorphisms to the terms of the resolution $S_\bullet(\frakg,\frakp)\otimes V$, 
one obtains an $U(\frakg)$-resolution 
$$
\displaylines{ 
\cdots \to U(\frakg)\otimes_{U(\frakp)} (\Lambda^2(\frakg/ \frakp) \otimes V\vert{}_\frakp) 
 \stackrel{d_2}{\to} U(\frakg)\otimes_{U(\frakp)} (\frakg/ \frakp \otimes V\vert{}_\frakp)  
\hfill\cr\hfill  
\stackrel{d_1}{\to} U(\frakg)\otimes_{U(\frakp)} V\vert{}_\frakp 
\stackrel{\varepsilon}{\to} V\to 0,
}  
$$
where the differentials $d_k$ are now given by
$$
\displaylines{ 
d_k(1\otimes \bar{x}_1\wedge\cdots\wedge \bar{x}_k\otimes v) 
:=  
\sum_{i=1}^k (-1)^{i-1}\, x_i\otimes 
\bar{x}_1\wedge\cdots\wedge \widehat{\bar{x}_i}\wedge \cdots \wedge \bar{x}_k \otimes v 
\hfill\cr\hfill
+ \sum_{1\leq i<j\leq k} (-1)^{i+j}\; 1\otimes \pi_\frakp([x_i,x_j])\wedge \bar{x}_1\wedge\cdots 
\wedge \widehat{\bar{x}_i}\wedge \cdots \widehat{\bar{x}_j}\wedge \cdots\wedge \bar{x}_k \otimes v 
\cr
\phantom{d_k(1\otimes \bar{x}_1\wedge\cdots\wedge \bar{x}_k\otimes v) :} 
+ \sum_{i=1}^k (-1)^i \; 1 \otimes 
\bar{x}_1\wedge\cdots\wedge \widehat{\bar{x}_i}\wedge \cdots \wedge \bar{x}_k \otimes x_i v,   
\hfill
}
$$
for $x_1,\dots,x_k\in\frakg$ and $v\in V$ (we have denoted $\pi_\frakp(x_i)$ by $\bar{x}_i$). 
We shall call it the standard resolution of $V$ relative to the pair $(U(\frakg), U(\frakp))$,  
and denote it by $S_\bullet(\frakg,\frakp,V)$. When $V = V_\ZZ(\lambda)$, we shall denote it 
by $S_\bullet(\frakg,\frakp,\lambda)$.

\subsection{}\label{weightsinwedge}
Let $p$ be a prime integer and recall the notation of \ref{alcove}. 

\begin{lem} Let $\lambda\in X^+\cap \barCp$. Then all weights $\nu$ of 
$V_\ZZ(\lambda)\otimes \Lambda(\frakg/\frakp)$ satisfy 
$\langle\nu+\rho,\alpha^\vee\rangle \leq p$, for all $\alpha\in R$.
\end{lem} 

\pf As $T$-module, $\Lambda(\frakg/\frakp)$ identifies with $\Lambda(\fraku_P^-)$ and 
hence is a submodule of $\Lambda(\fraku^-)$, where $\fraku^-$ is the Lie algebra of 
$U^-$. 

By a result of Kostant (\cite[Lemma 5.9]{Ko1}), there is a $T$-isomorphism 
$$
\rho\otimes \Lambda(\fraku^-) \cong V_\ZZ(\rho).
$$
 Therefore, if $\nu$ is a weight of $V_\ZZ(\lambda)\otimes \Lambda(\frakg/\frakp)$, then 
$\nu+\rho$ is a weight of $V_\ZZ(\lambda)\otimes V_\ZZ(\rho)$. 
This implies that 
$\langle\nu+\rho,\alpha^\vee\rangle \leq p$, for all $\alpha\in R$. 

Indeed, let 
$\mu$ be the dominant $W$-conjugate of $\nu+\rho$, it is also a weight of 
$V_\ZZ(\lambda)\otimes V_\ZZ(\rho)$. Clearly, it suffices to prove that 
$\langle\mu,\alpha^\vee\rangle \leq p$, for all $\alpha\in R^+$. 
Further, since $\mu$ is dominant, it suffices to prove that 
$\langle\mu,\gamma^\vee\rangle \leq p$ when $\gamma^\vee$ is a maximal coroot. 
But it is well-known that a maximal coroot is a dominant coweight, i.e. satisfies 
$\langle\beta,\gamma^\vee\rangle\geq 0$ for all $\beta\in R^+$, see e.g. 
\cite[VI,\S\,1, Prop.8]{Bki}. Finally, since $\mu = \lambda+\rho-\theta$ with 
$\theta\in \NN R^+$, it follows that 
$$
\langle\mu,\gamma^\vee\rangle \leq \langle\lambda+\rho,\gamma^\vee\rangle\leq p . 
$$
This proves the lemma.

\subsection{Verma modules and filtrations.}\label{defVerma}\label{Vermafilt}  
For any $\xi\in X_L^+$, define the generalized Verma module (for $U(\frakg)$ and $U(\frakp)$) 
$$
M_\frakp(\xi) := U(\frakg)\otimes_{U(\frakp)} V_\ZZ^L(\xi).
$$
For any commutative ring $A$, set $M_\frakp^A(\xi) := M_\frakp(\xi)\otimes A$ and observe that 
it identifies with  
$U_A(\frakg)\otimes_{U_A(\frakp)} V_A^L(\xi)$. 

For $\lambda\in X^+$, we set also 
$$
S_\bullet^A(\frakg,\frakp,\lambda) := S_\bullet(\frakg,\frakp,\lambda)\otimes A.
$$

Let $i\geq 0$. It follows from Lemmas \ref{weightsinwedge} and \ref{donkin} that 
there exists a $P$-module filtration 
$$
0 = F_0 \subset \cdots \subset F_r = \Lambda^i(\frakg/\frakp)\otimes V_\Zp(\lambda)
$$ 
such that each $F_j/F_{j-1}$ is isomorphic to $V_\Zp^L(\xi_j^i)$, for some $\xi_j^i\in X_L^+$ 
(not necessarily distinct). Let us denote by $\Omega_\frakp^i(\lambda)$ the 
multiset of those $\xi_j^i$ (each $\xi\in X_L^+$ occuring as many times as $V_\Zp^L(\xi)$ 
occurs in the filtration). 

Moreover, as $U(\frakg)$ is free over $U(\frakp)$, the functor $U(\frakg)\otimes_{U(\frakp)}-$ 
is exact. Therefore, one obtains the 

\begin{cor} Let $\lambda\in X^+\cap \barCp$. Then each $S_i^\Zp(\frakg,\frakp,\lambda)$ 
admits a finite filtration by $U_\Zp(\frakg)$-modules such that the successive quotients are 
the $M_\frakp^\Zp(\xi)$, for $\xi\in \Omega_\frakp^i(\lambda)$. 
\end{cor}

\subsection{The Carter-Lusztig algebra.}\label{Carter-Lusztig} 
In this subsection, we assume that $G = GL_n$. 
Let $\{E_{ij}\}_{i,j = 1}^n$ be the standard basis (over $\ZZ$) of $\frakg = \frakgl_n$ 
and let $c$ denote the central element $\sum_i E_{ii}$. Let 
$U'_\ZZ(\frakgl_n)$ denote the subalgebra of $U_\QQ(\frakgl_n)$ generated by $\frakgl_n$ 
and by the elements 
$$
\binom{c}{r} := \frac{1}{r!} \prod_{i=1}^r (c-i+1),  
$$
for $r\geq 1$. Clearly, these elements are invariant for the adjoint action of $GL_n$. 
We call $U'_\ZZ(\frakgl_n)$  the Carter-Lusztig algebra (see \cite{CL}).

Recall that the dominant weights for $GL_n$ identify with 
sequences 
of integers $\lambda_1\geq \cdots \geq \lambda_n$; the corresponding Weyl module $V_\ZZ(\lambda)$  
identifies with the submodule generated by a highest weight vector in the tensor product 
$$
\bigotimes_{i=1}^{n-1} (\Lambda^i V_\ZZ)^{\otimes (\lambda_i-\lambda_{i+1})} \bigotimes 
\det{}^{\otimes \lambda_n}, 
$$
where $V_\ZZ$ denotes the natural representation of $GL_n$. 

Set $\vert\lambda\vert := \sum_{i=1}^n \lambda_i$. Then the center $Z \cong \GG_m$ of $GL_n$ 
acts on $V_\ZZ(\lambda)$ with the weight $\vert\lambda\vert$ and hence  
each $\binom{c}{r}$ acts on $V_\ZZ(\lambda)$ by the 
integer $\binom{\vert\lambda\vert}{r}$ (see \cite[3.8, Remark 2]{CL}). More generally, 
as $Z$ acts trivially on the $P$-module $\Lambda^\bullet(\frakg/\frakp)$, it also acts on 
$\Lambda^\bullet(\frakg/\frakp)\otimes V_\ZZ(\lambda)$ by the weight $\vert\lambda\vert$. 

Also, let $U'_\ZZ(\frakp)$ (resp. $U'_\ZZ(\frakt)$) be the subalgebra of $U_\QQ(\frakg)$ 
generated by $\frakp$ (resp. $\frakt$) and the $\binom{c}{r}$, $r\geq 1$. Then, one deduces 
from the PBW theorem that 
$$
U'_\ZZ(\frakg) \cong U_\ZZ(\fraku_P^-) \otimes U'_\ZZ(\frakp).
$$

Let $\lambda\in X^+$. Then $\Lambda^\bullet(\frakg/\frakp)\otimes V_\ZZ(\lambda)$ is 
an $\calU_\ZZ(P)$-module and hence an $U'_\ZZ(\frakp)$-module. Clearly, there is 
an isomorphism of $U_\ZZ(\frakg)$-modules 
$$
U_\ZZ(\frakg)\otimes_{U_\ZZ(\frakp)} 
\bigl(\Lambda^\bullet(\frakg/\frakp)\otimes V_\ZZ(\lambda) \bigr) 
\cong 
U'_\ZZ(\frakg)\otimes_{U'_\ZZ(\frakp)} 
\bigl(\Lambda^\bullet(\frakg/\frakp)\otimes V_\ZZ(\lambda) \bigr). 
$$
It follows that $S_\bullet(\frakgl_n,\frakp,\lambda)$ carries a natural action of 
$U'_\ZZ(\frakgl_n)$, where each $\binom{c}{r}$ acts by the integer $\binom{\vert\lambda\vert}{r}$. 

In the sequel, we shall have to work either with the algebra $U'_\ZZ(\frakgl_n)$, when 
$G = GL_n$, or with the algebra $U_\ZZ(\frakg)$, when $G$ is simply connected and not of type $A$. 
In order to have uniform notation, we shall set in the latter case $U'_\ZZ(\frakg) = U_\ZZ(\frakg)$ 
and $U'_\ZZ(\frakt) = U_\ZZ(\frakt)$.

\subsection{Central characters.}\label{centralchar} 
Let $A$ be a commutative ring. Denote by $\fraku^-$ and $\fraku$ the Lie algebras 
of $U^-$ and $U$, respectively. By the PBW theorem, one has 
$$
U'_A(\frakg) = U'_A(\frakt) \oplus 
\left(\fraku^- U'_A(\frakg) + U'_A(\frakg)\fraku \right). 
$$
Let $\delta_A$ denote the $A$-linear 
projection from $U'_A(\frakg)$ to $U'_A(\frakt)$ defined by this decomposition.

Let $U'_A(\frakg)^G \subset U'_A(\frakg)^T$ be the subrings of $G$-invariant and $T$-invariant 
elements for the adjoint action. 
Denote by $\theta_A$ the restriction of $\delta_A$ to 
$U'_A(\frakg)^T$. Then $\theta_A$ is a ring homomorphism; indeed one sees easily that the arguments 
in the proof of \cite[Lemme 7.4.2]{Di} or \cite[Lemma 5.1]{KW} carry over in our case. 

For any $\mu\in X(T)$, its differential $d\mu$ induces an $A$-linear map $\frakt_A\to A$ and 
hence an $A$-algebra morphism $U_A(\frakt) \to A$, still denoted by $d\mu$. Further, in the case 
where $G = GL_n$, one sees easily that $d\mu$ extends to an $A$-algebra morphism 
$U'_A(\frakt) \to A$ : its value on an element $\binom{c}{r}$ is the image in $A$ of the integer 
$\binom{\vert\mu\vert}{r}$. 
Thus, $\mu$ gives 
rise to an $A$-algebra morphism $\chi_{\mu,A} := d\mu\circ \theta_A$, 
from $U'_A(\frakg)^G$ to $A$. 

For any morphism of commutative rings $f : A\to B$, it is easily seen that the following diagram 
is commutative:
$$
\begin{array}{ccccc}
U_A(\frakg)^G & \stackrel{\theta_A}{\to} & U'_A(\frakt) & \stackrel{d\mu}{\to} & A \\
\phantom{f} \downarrow f & & \phantom{f} \downarrow f & & \phantom{f} \downarrow f \\
U_B(\frakg)^G & \stackrel{\theta_B}{\to} & U'_B(\frakt) & \stackrel{d\mu}{\to} & B .
\end{array} 
$$
Thus, one has $\chi_{\mu,B} = f\circ \chi_{\mu,A}$. 

\medskip {\it Remark}. The map $U'_A(\frakg)^G \otimes B \to U'_B(\frakg)^G$ need not be surjective 
(for example, when $G=SL_n$, $A = \ZZ$ and $B= \Fp$ where $p$ divides $n$). 

\smallskip  
Since elements of $U'_A(\frakg)^T$ have weight zero, one has 
$$
U'_A(\frakg)^T \subseteq U'_A(\frakt)\oplus \fraku^- U'_A(\frakg) \fraku . 
$$
Therefore, if $M$ is an $U'_A(\frakg)$-module generated by an element $v$ of weight 
$\mu$ annihilated by $\fraku$, then $U'_A(\frakg)^G$ acts on $M$ by the character 
$\chi_{\mu,A}$ (see \cite[Prop.7.4.4]{Di}).

Let $\pi$ denote the morphism $\Zp\to \barFp$ and let $\chi_{\mu,p} := \chi_{\mu,\Zp}$ 
and $\barchi_{\mu,p} := \pi\circ \chi_{\mu,p} = \chi_{\mu,\barFp}$. Set also 
$J_{\mu,p} := \Ker\chi_{\mu,p}$.

\subsection{Decomposition w.r.t. central characters mod. $p$.}\label{blocks}  

Let $\lambda\in X^+$ and let $p$ be a prime integer such that $\lambda\in \barCp$. 
Recall the multisets $\Omega_\frakp^i(\lambda)$ from \ref{Vermafilt} and let 
$\Omega_\frakp^\bullet(\lambda)$ denote their disjoint union.  

By Corollary \ref{Vermafilt}, each $S_i^\Zp(\frakg,\frakp,\lambda)$ admits a finite 
$U_\Zp(\frakg)$-filtration, whose quotients are the $M_\frakp^\Zp(\xi)$, where $\xi$ 
runs through $\Omega_\frakp^i(\lambda)$. It follows that $S_\bullet^\Zp(\frakg,\frakp,\lambda)$ 
is annihilated by the ideal 
$$
 I := \prod_{\xi\in \Omega_\frakp^\bullet(\lambda)} J_{\xi,p}
$$
(each $\xi$ being counted with its multiplicity). 

\smallskip The following lemma is straightforward. 

\begin{lem} Let $A$ be a commutative ring and $P_1,\dots,P_r$ ideals of $A$ such that 
$P_1\cdots P_r = 0$ and $P_i+P_j = A$ if $j\not= i$. Then, for any $A$-module $M$, 
one has 
$$
M = \bigoplus_{i=1}^r M^{P_i}, \quad\mbox{where } M^{P_i} = \{m\in M\mid P_i\, m=0\}.
$$
Further, the assignment $M\mapsto M^{P_i}$ is an exact functor. 
\end{lem}

We shall apply the lemma to $A := U'_\Zp(\frakg)^G/I$. 
Note that $A$ is a finite $\Zp$-module. Moreover, it is easily seen that 
the maximal ideals of $A$ are the $pA + J_{\xi,p} = \Ker \barchi_{\xi,p}$. (By abuse of 
notation, we still denote by $J_{\xi,p}$ the image of $J_{\xi,p}$ in $A$). Let  
$\barchi_1,\dots,\barchi_r$ be the distinct algebra homomorphisms $A\to \barFp$, with 
$\barchi_1 = \barchi_{\lambda,p}$ and, for $i = 1,\dots,r$, let 
$$
P_i := \prod_{{\xi\in\Omega_\frakp^\bullet(\lambda) \atop \barchi_{\xi,p} = \barchi_i}} J_{\xi,p} .
$$
Clearly, $P_1\cdots P_r = 0$ and $pA + P_i + P_j = A$ if $j\not= i$. Since $A$ is a finite 
$\Zp$-module, the latter implies, by Nakayama lemma, that $P_i + P_j = A$ if $j\not= i$. 

Then, one deduces from the previous lemma that each $S_i^\Zp(\frakg,\frakp,\lambda)$ 
contains as a direct summand the $U_\Zp(\frakg)$-submodule 
$$
S_i^\Zp(\frakg,\frakp,\lambda)_{\barchi_{\lambda,p}} := S_i^\Zp(\frakg,\frakp,\lambda)^{P_1}.
$$ 
Moreover, since the differentials in the complex $S_\bullet^\Zp(\frakg,\frakp,\lambda)$ are 
$U_\Zp(\frakg)$-equivariant, $S_\bullet^\Zp(\frakg,\frakp,\lambda)_{\barchi_{\lambda,p}}$ 
is a direct summand {\it subcomplex}. 

Further, since $M\mapsto M_{\barchi_{\lambda,p}}$ is an exact functor and since 
$$
M_\frakp^\Zp(\xi)_{\barchi_{\lambda,p}} = \left\{ 
\begin{array}{cl} 
M_\frakp^\Zp(\xi) & \mbox{if } \barchi_{\xi,p} = \barchi_{\lambda,p}\,;\\
0 & \mbox{otherwise},
\end{array}
\right.
$$
one obtains, as in \cite[Lemma 9.7]{BGG}, the following 

\begin{cor} $S_\bullet^\Zp(\frakg,\frakp,\lambda)$ contains 
the subcomplex 
$S_\bullet^\Zp(\frakg,\frakp,\lambda)_{\barchi_{\lambda,p}}$ 
as a direct summand. Moreover, 
for $i\geq 0$ each 
$S_i^\Zp(\frakg,\frakp,\lambda)_{\barchi_{\lambda,p}}$ has a filtration 
whose quotients are the $M_\frakp^\Zp(\xi)$, for those 
$\xi\in \Omega_\frakp^i(\lambda)$ $($counted with multiplicities$)$ such that 
$\barchi_{\xi,p} = \barchi_{\lambda,p}$. 
\end{cor}

\subsection{}\label{weakBGG} 
 
The first step towards the description of 
$S_\bullet^\Zp(\frakg,\frakp,\lambda)_{\barchi_{\lambda,p}}$ 
is the following proposition. 

\begin{prop} Let $\xi\in\Omega_\frakp^\bullet(\lambda)$. Suppose that 
$\barchi_{\xi,p} = \barchi_{\lambda,p}$. Then $\xi = w\cdot \lambda$ for some 
$w\in W^L$. 
\end{prop}  

\pf Let $\xi$ be as in the proposition. 
Consider the following two cases. 

1) If $G = GL_n$, the fact that $\barchi_{\xi,p} = \barchi_{\lambda,p}$ implies that 
$\xi\in W_p\cdot \lambda$, by the proofs of Theorems 3.8 and 4.1 in \cite{CL}. 

2) If $G$ is quasi-simple and simply-connected, $\barchi_{\xi,p} = \barchi_{\lambda,p}$ implies, 
by \cite[Th.\,2]{KW}, that there exist $y\in W$ and $\nu\in X(T)$ such that 
$y\cdot \xi = \lambda + p\nu$. Moreover, 
since $y\cdot\xi$ is a weight of $\Lambda(\frakg/\frakp)\otimes V_\ZZ(\lambda)$, then 
$y\cdot\xi - \lambda\in Q(R)$ and hence $p\nu\in Q(R)\cap pX(T)$. 
But, when $R$ is irreducible and not of type $A$, 
our assumption that 
$$
p\geq \langle\lambda+\rho,\gamma^\vee\rangle \geq \langle\rho,\gamma^\vee\rangle, 
$$
where $\gamma^\vee$ is the highest coroot, implies that $p > \vert X(T)/Q(R)\vert $. 
It follows that $\nu\in Q(R)$ and hence $\xi\in W_p\cdot\lambda$. 

Thus, in both cases, we have obtained that $\xi\in W_p\cdot\lambda$. Now, let 
$w\in W$ such that $w^{-1}(\xi+\rho)$ is dominant and let 
$\xi^+ := w^{-1}\cdot\xi$. Then, by Lemma \ref{weightsinwedge}, 
$\xi^+\in \barCp$. But $\xi^+\in W_p\cdot\lambda$; since $\barCp$ is a fundamental domain for the 
dot action of $W_p$, it follows that $\xi^+ = \lambda$, and hence $\xi = w\cdot \lambda$. 

Further, since $\xi\in\Omega_\frakp^\bullet(\lambda)\subseteq X_L^+$, for any $\alpha\in R_L^+$ 
one has $\langle w\cdot\lambda, \alpha^\vee \rangle \geq 0$ and hence 
$$
\langle \lambda+\rho, w^{-1}\alpha^\vee \rangle \geq \langle \rho, \alpha^\vee \rangle > 0.
$$
This implies that $w\in W^L$. The proposition is proved.

\medskip 
We can now prove the following analogue of \cite[Th.\,9.9]{BGG} and \cite[Th.\,3.10]{Le}, 
\cite[Th.\,7.11]{Rocha}. 

\begin{th} Suppose that $\lambda\in X^+\cap \barCp$. 
Then 
$S_\bullet^\Zp(\frakg,\frakp,\lambda)_{\barchi_{\lambda,p}}$ is an 
$U'_\Zp(\frakg)$-resolution of $V_\Zp(\lambda)$ and each 
$S_i^\Zp(\frakg,\frakp,\lambda)_{\barchi_{\lambda,p}}$ with $i\geq 0$ has a filtration 
whose quotients are exactly the $M_\frakp^\Zp(w\cdot\lambda)$, for $w\in W^L(i)$, 
each occuring once. 
\end{th}

\pf 
By Corollary \ref{blocks}, combined with the previous proposition, each 
$S_i^\Zp(\frakg,\frakp,\lambda)_{\barchi_{\lambda,p}}$ with $i\geq 0$ has a filtration 
whose quotients are the $M_\frakp^\Zp(\xi)$, for those 
$\xi\in \Omega_\frakp^i(\lambda)$ $($counted with multiplicities$)$ such that 
$\xi = w\cdot\lambda$ for some $w\in W^L$. 

Conversely, for $w\in W^L$, Kostant has showed that $V_\QQ^L(w\cdot \lambda)$ occurs with 
multiplicity one in $\Lambda^\bullet(\frakg/\frakp)\otimes V_\ZZ(\lambda)$, 
in degree equal to $\ell(w)$, see [Ko1], Lemma 5.12 and end of proof of 
Th. 5.14. This completes the proof of the theorem.

\subsection{Proof of theorem \ref{kostant}.}\label{pf-kostant} 

As $U_\Zp(\fraku_P^-)$-module, each $M_\frakp^\Zp(\xi)$ is isomorphic to 
$U_\Zp(\fraku_P^-)  \otimes V_\Zp^L(\xi)$, hence free. Thus, 
by the previous theorem, 
$S_i^\Zp(\frakg,\frakp,\lambda)_{\barchi_{\lambda,p}}$ is a free $U_\Zp(\fraku_P^-)$-module, 
for each $i\geq 0$. 

Therefore, $H_\bullet(\fraku_P^-, V_\Zp(\lambda))$ is the homology of the complex 
$$
C_\bullet := 
\Zp\otimes_{U_\Zp(\fraku_P^-)} S_\bullet^\Zp(\frakg,\frakp,\lambda)_{\barchi_{\lambda,p}}.  
$$
Further, by the previous theorem, again, for $i\geq 0$ each $C_i$ 
has an $L$-module filtration whose successive quotients are the 
$V_\Zp^L(w\cdot\lambda)$, for $w\in W^L(i)$. 

By Corollary \ref{vanish}, applied to $L$, one obtains that these filtrations split, that is, 
for each 
$i\geq 0$ one has isomorphisms of $L$-modules 
$$
C_i \cong \oplus_{w\in W^L(i)} V_\Zp^L(w\cdot\lambda). 
$$

Further, we claim that the differentials $d_i : C_i\to C_{i-1}$ are zero. Indeed, one has 
$H_i(C_\bullet)\otimes \QQ \cong H_i(\fraku_P^-, V_\QQ(\lambda))$ and, by Kostant's theorem 
(\cite[Cor\,8.1]{Ko1} or \cite{BGG}, Cor. of Th.\,9.9), 
the latter is isomorphic to $C_i\otimes \QQ$. It follows, 
for a reason of dimension, that $d_i\otimes 1 = 0$. Since $C_{i-1}$ is torsion-free, 
this implies that $d_i = 0$.  
 
Thus, we have obtained, for each $i\geq 0$, an isomorphism of $L$-modules 
$$
H_i(\fraku_P^-, V_\Zp(\lambda)) \cong \bigoplus_{w\in W^L(i)} V_\Zp^L(w\cdot\lambda).
$$
This completes the proof of Theorem \ref{kostant}.

\subsection{Analogue in cohomology.}\label{cohom-u}
Recall the anti-involution $\tau$ from \ref{def-tau}; it exchanges 
$P^-$ and $P$ and stabilizes $L$. Let $\lambda\in X^+\cap\barCp$. 
Since $H_\bullet(\fraku_P^-, V)$ is a free $\Zp$-module, one obtains, by 
standard arguments, an isomorphism of $L$-modules  
$$
H_\bullet(\fraku_P^-, V_\Zp(\lambda))^\tau \cong H^\bullet(\fraku_P, V_\Zp(\lambda)^\tau).
$$ 
Further, since $V_\Zp(\lambda) = V_\Zp(\lambda)^\tau$ and 
$V_\Zp^L(w\cdot\lambda) = V_\Zp^L(w\cdot\lambda)$, for $w\in W^L$,  
by Corollary \ref{irred}, applied to $G$ and $L$, one obtains the 

\begin{cor} Let $\lambda\in X^+\cap\barCp$. For each $i\geq 0$, there is 
an isomorphism of $L_\Zp$-modules 
$$
H^i(\fraku_P, V_\Zp(\lambda)) \cong \bigoplus_{w\in W^L(i)} V_\Zp^L(w\cdot\lambda).
$$
\end{cor}

\section{Cohomology of the groups $U_P^-({\bf Z})$} 

\subsection{}\label{N-series}\label{Hartley}  
Let $\Gamma$ be a finitely generated, torsion free, nilpotent group.    
Let $\calF$ be a decreasing sequence 
$\Gamma = F^1\Gamma \supseteq F^2\Gamma \supseteq \cdots$ of normal subgroups 
of $\Gamma$. 
Following the terminology in Passman's book \cite[p.85]{Pas}, let us say that $\calF$  
is an $N$-series if $(F^i\Gamma, F^j\Gamma)\subseteq F^{i+j}\Gamma$ for all $i,j$. 
Further, $\calF$ is called an $N_0$-series if it is an $N$-series and 
each $F^i\Gamma/F^{i+1}\Gamma$ is torsion-free.

If $\calF$ is an $N$-series, the associated graded abelian group 
$$
\gr_\calF \Gamma := \bigoplus_{i\geq 1} F^i\Gamma/F^{i+1}\Gamma
$$
has a natural structure of Lie algebra over $\ZZ$ (see, for example, 
\cite[Chap.\,I, Th.\,2.1]{Laz}). 
 By \cite[Th.\,2.2]{Jen}, if $\calF$ is an $N_0$-series, $\gr_\calF \Gamma$ 
is a finite free $\ZZ$-module; its rank, $r$, is an invariant called 
the rank of $\Gamma$.

\smallskip 
For $i\geq 1$, let $\{C^i(\Gamma)\}$ denote the lower central 
series; as is well-known, it is the fastest descending $N$-series. 
We shall denote the corresponding graded Lie algebra simply by $\gr \Gamma$. 
Further, set 
$$
C^{(i)}(\Gamma) := \{x\in \Gamma \mid x^n \in C^i(\Gamma) \mbox{ for some }n>0\}.
$$
By \cite[Chap.\,11, Lemma 1.8]{Pas} (see also \cite[\S\,4]{Hall}), $\{C^{(i)}(\Gamma)\}$ 
is an $N_0$-series. It is clearly the fastest descending $N_0$-series. 
Following \cite[\S\,4]{Hall}, we will call it the isolated lower central series.  
We will denote by 
$\grisol\Gamma$ the associated Lie algebra over $\ZZ$ 
$$
\grisol\Gamma := \bigoplus_{i\geq 1} C^{(i)}(\Gamma)/C^{(i+1)}(\Gamma),
$$
which is also a free $\ZZ$-module of rank $r$. Clearly, there is an isomorphism of 
graded Lie algebras $\gr \Gamma \otimes \QQ \cong \grisol\Gamma \otimes \QQ$. 

\medskip 
Let $I$ denote the augmentation ideal of the group ring $\ZZ\Gamma$ and, for $n\geq 0$, 
let $I^{(n)}$ denote the isolator of $I^n$, that is, 
$$
I^{(n)} := \{x\in\ZZ\Gamma \mid mx\in I^n \mbox{ for some } m > 0\}.
$$
Equivalently, if $I_\QQ$ denotes the augmentation ideal of $\QQ\Gamma$, then 
$I^{(n)} = \ZZ\Gamma\cap I_\QQ^n$. 

Let us consider the graded rings 
$$
\grisol\ZZ\Gamma := \bigoplus_{n\geq 0} I^{(n)}/I^{(n+1)} 
\quad \mbox{and}\quad 
\gr\QQ\Gamma := \bigoplus_{n\geq 0} I_\QQ^n/I_\QQ^{n+1}.
$$
The former is a subring of the latter and, by a result of Quillen (\cite{Qui}), 
there is an isomorphism of graded Hopf algebras 
$U_\QQ(\gr\Gamma \otimes \QQ) \cong \gr\QQ\Gamma$. Further, one has the following 
more precise result of Hartley :

\begin{th}{\rm (\cite[Th.\,2.3.3$'$]{Hart2})}  There is 
an isomorphism of graded Hopf algebras 
$$U_\ZZ(\grisol\Gamma) \cong \grisol\ZZ\Gamma.$$ 
\end{th}

\subsection{}\label{grisol-u}\label{Vergne}

Let $A$ be a finitely generated subring of $\QQ$ (thus, $A = \ZZ[1/m]$ for some 
$m$ and $A$ is a PID). Let $\fraku$ be a nilpotent Lie algebra over $A$, which 
is a finite free $A$-module, say of rank $r$. Let $\fraku_\QQ = \fraku\otimes_A \QQ$, 
then $U_\QQ(\fraku_\QQ)\cong U_A(\fraku)\otimes_A \QQ$; we shall denote it by 
$U_\QQ(\fraku)$. 
By the PBW theorem, 
$U_A(\fraku)$ is a subalgebra of $U_\QQ(\fraku)$. 

\smallskip 
Let $\calF$ be a decreasing sequence 
$\fraku = F^1\fraku \supseteq F^2\fraku \supseteq \cdots$ of Lie ideals 
of $\fraku$. As in the previous paragraph, let us say that $\calF$ is an $N$-series 
if $[F^i\fraku, F^j\fraku]\subseteq F^{i+j}\fraku$, and is an $N_0$-series if further 
each $F^i\fraku / F^{i+1}\fraku$ (which is a finitely generated module over the PID $A$) 
is torsion free, and hence a free $A$-module. 

Let $C^i(\fraku)$ denote the lower central series of $\fraku$ and 
define the isolated lower central series 
$\{C^{(i)}(\fraku)\}$ by 
$$
C^{(i)}(\fraku) := \{x\in \fraku \mid nx\in C^i(\fraku) \mbox{ for some }n>0\}.
$$
This is, clearly, the fastest descending $N_0$-series of $\fraku$. 
Consider the graded Lie algebras  
$$
\grisol \fraku := \bigoplus_{i\geq 1} C^{(i)}(\fraku)/C^{(i+1)}(\fraku)
\quad \mbox{and}\quad 
\gr \fraku_\QQ := \bigoplus_{i\geq 1} C^{i}(\fraku_\QQ)/C^{i+1}(\fraku_\QQ). 
$$
Then $\grisol \fraku$ is a free $A$-module of rank $r$ and there is an isomorphism 
of graded Lie algebras $(\grisol \fraku)\otimes_A \QQ \cong \gr \fraku_\QQ$.   

\smallskip Let $J_\QQ$ denote the augmentation ideal of $U_\QQ(\fraku)$. Then the 
graded algebra 
$$
\gr U_\QQ(\fraku) := \bigoplus_{n\geq 0} J_\QQ^n/J_\QQ^{n+1}
$$
is a primitively generated, graded Hopf algebra; it is isomorphic to 
$U_\QQ(\gr \fraku_\QQ)$, by \cite{Knus} or \cite[Prop.\,1]{Ver}. 
In fact, as in the case of group rings, a little more is true. 
For $n\geq 1$, let 
$J^{(n)} = U_A(\fraku)\cap J_\QQ^n$. Then the graded ring 
$$
\grisol U_A(\fraku) := \bigoplus_{n\geq 0} J^{(n)}/J^{(n+1)}
$$
identifies with a subring of $\gr U_\QQ(\fraku)$. Further, one deduces from the 
proof of \cite[Prop.\,1]{Ver} the following result. 
Let $X_1,\dots,X_r$ be an $A$-basis of $\fraku$ compatible with the filtration 
$\{C^{(i)}(\fraku)\}$, {\it i.e.}, such that for $s = 1,\dots,c$, the 
$X_j$ with $j > r - \dim C^s(\fraku_\QQ)$ form an $A$-basis of $C^{(s)}(\fraku)$, 
 and, for each $i$, let $\mu(i)$ be the largest integer $k$ 
such that $X_i\in C^{(k)}(\fraku)$.

\begin{prop} 

{\rm a)} For any $n\geq 0$, the ordered monomials 
$X_1^{n_1}\cdots X_r^{n_r}$ with $\sum_{i=1}^r n_i\,\mu(i) \geq n$ 
form an $A$-basis of $J^{(n)}$. 

\vskip3pt \noindent {\rm b)} 
There is an isomorphism of graded Hopf algebras 
$U_A(\grisol \fraku) \cong \grisol U_A(\fraku)$. 
\end{prop}

\subsection{}\label{same-gr}\label{Hall}

Let $\Gamma = H_1 \supset \cdots \supset H_{r+1} = \{1\} $ be a refinement of the 
isolated lower central series such that each $H_i/H_{i+1}$ is an infinite cyclic group, 
generated by the image of an element $g_i$ of $H_i$. Then,  
$\{g_1,\dots,g_r\}$ is called a system of canonical parameters (or canonical basis) 
of $\Gamma$; it induces a bijection $\ZZ^r \cong \Gamma$, given by 
$(n_1,\dots,n_r)\mapsto g_1^{n_1}\cdots g_r^{n_r}$; we will denote the R.H.S. simply by 
$g(n_1,\dots,n_r)$. Let $\{e_1,\dots,e_r\}$ be the standard basis of $\ZZ^r$, then 
$g(e_i) = g_i$.  

\medskip 
Let $\calP_{r,r}$ denote the subring of the polynomial ring 
$\QQ[\xi_1,\dots,\xi_r,\eta_1,\dots, \eta_r]$ consisting of 
those polynomials which take integral values on $\ZZ^r\times \ZZ^r$. 
By a result of Ph. Hall \cite[Th.\,6.5]{Hall}, there exist polynomials 
$P_1,\dots, P_r\in \calP_{r,r}$ such that 
$$
g(x_1,\dots,x_r)\, g(y_1,\dots,y_r) = g(P_1(x,y),\dots,P_r(x,y)), \leqno (\star)
$$
for any $x,y\in\ZZ^r$. 

Therefore, there exists an algebraic unipotent group 
scheme $U$, defined over a finitely generated subring $A$ of the rationals, 
and whose underlying scheme is affine space ${\Bbb A}_A^r$, such that 
$\Gamma$ identifies with the subgroup $\ZZ^r$ of $U(A) = A^r$. 

\medskip\noindent {\it Remark}. If $\Gamma$ is of class $c$, 
one may take $A = \ZZ[1/c!]$; this can be deduced, for example, from 
the Campbell-Hausdorff formula. 

\medskip 
Let $k \in \{1,\dots,r\}$. Since $P_k(x,0) = x_k$ and $P_k(0,y) = y_k$ for every 
$x,y\in \ZZ^r$, the part of degree $\leq 1$ of $P_k$ is $\xi_k + \eta_k$ and its part of degree  
$2$, call it $b_k$, is bilinear in the $\xi_i$ and the $\eta_j$. Thus, one has
$$
P_k(\xi,\eta) = \xi_k+\eta_k + \sum_{i,j=1}^r b_k(e_i,e_j)\, \xi_i \eta_j + \mbox{terms of degree $> 2$}.
$$ 

 Let $\frakm$ denote the ideal $(\xi_1,\dots,\xi_r)$ of 
$A[U] = A[\xi_1,\dots,\xi_r]$, let 
$$
\fraku := \Hom_A(\frakm/\frakm^2, A)
$$ 
be the Lie algebra of $U$ over $A$, and let $\{v_1,\dots, v_r\}$ be the  
$A$-basis of $\fraku$ dual to the basis $\{\bar{\xi}_1,\dots,\bar{\xi}_r\}$. 
Then, the Lie brackets are given by 
$$
[v_i,v_j] = \sum_{k=1}^r (b_k(e_i,e_j) - b_k(e_j,e_i))\, v_k, \leqno (1)
$$
see, for example, \cite[\S\,1]{LP} or \cite[\S\,1]{CP}.

\begin{prop} There is an isomorphism of graded Lie algebras over $A$ 
$$
\grisol \Gamma \otimes_\ZZ A \cong \grisol \fraku, 
$$
under which each $\bar{g}_i$ corresponds to $\bar{v}_i$. 
\end{prop}

\pf First, for each $i$, let $\nu(i)$ denote the largest integer $n$ 
such that $g_i \in C^{(n)}(\Gamma)$. 
Denote by $\bar{g}_i$ the image of 
$g_i$ in $\grisol^{\nu(i)}\Gamma$; then $\{\bar{g}_1,\dots,\bar{g}_r\}$ is a 
$\ZZ$-basis of $\grisol\Gamma$. 

For $k = 1,\dots,r$, let $Q_k := P_k - \xi_k - \eta_k$ be the part of $P_k$ of degree $>1$.  
Recall that, 
for $x_1,\dots,x_r\in \ZZ$,  $g(\sum_{i=1}^r x_i\,e_i)$ denotes the element 
$g_1^{x_1}\cdots g_r^{x_r}$ of $\Gamma$.

Let $i,j\in\{1,\dots,r\}$ be arbitrary with $i< j$. 
Then, for every $x,y\in \ZZ^r$, one has 
$g(xe_i)g(ye_j) = g(xe_i + ye_j)$ and hence $Q_k(xe_i,ye_j) = 0 = 
b_k(xe_i,ye_j)$ for any $k$. In particular, $b_k(e_i,e_j) = 0$. 

On the other hand, since 
$g_j^x \in C^{(\nu(j))}(\Gamma)$ and $g_i^y\in C^{(\nu(i))}(\Gamma)$ 
one has,  
$$
g_j^x\, g_i^y \equiv g_i^y\, g_j^x\, 
g\left(\sum_{{k\atop \nu(k)=\nu(i)+\nu(j)}} Q_k(x,y)e_k \right) 
\quad \mbox{mod. } C^{(\nu(i)+\nu(j)+1)}(\Gamma).
$$
Further, since the commutator induces a bilinear map on  
$\grisol \Gamma$, one has, when $\nu(k) = \nu(i)+\nu(j)$,  
$$
Q_k(xe_j,ye_i) = x y Q_k(e_j,e_i) = x y b_k(e_j,e_i). 
$$
Then, an easy computation shows that 
$$
g_i^x\,g_j^y\,g_i^{-x}\,g_j^{-y} \equiv 
g\left(\sum_{{k\atop \nu(k)=\nu(i)+\nu(j)}} -xy b_k(e_j,e_i)\,e_k \right) 
\quad \mbox{mod. } C^{(\nu(i)+\nu(j)+1)}(\Gamma).
$$
Using the fact that $b_k(e_i,e_j) = 0$, one 
deduces that the Lie bracket on $\grisol\Gamma$ is given by 
$$
[\bar{g}_i, \bar{g}_j] = 
\sum_{{k \atop \nu(k) = \nu(i)+\nu(j)}} (b_k(e_i,e_j)-b_k(e_j,e_i)) 
\, \bar{g}_k. 
\leqno (2)
$$

The proposition is then a consequence of the following claim. 

\smallskip\noindent {\bf Claim.} For each $\ell$, $C^{(\ell)}(\fraku)$ is the $A$-span of those 
$v_k$ such that $\nu(k)\geq \ell$. 

\smallskip 
Indeed, using $(1)$, the claim implies that $\grisol\fraku$ is the Lie algebra 
having an $A$-basis $\{\bar{v}_1,\dots, \bar{v}_r\}$ and 
brackets given by 
$$
[\bar{v}_i,\bar{v}_j] = \sum_{{k \atop \nu(k) = \nu(i)+\nu(j)}} 
(b_k(e_i,e_j) - b_k(e_j,e_i))\, \bar{v}_k .\leqno (3) 
$$ 
Comparing with $(2)$, one obtains that 
$\grisol \Gamma \otimes_\ZZ A \cong \grisol \fraku$. 

\smallskip 
Let us now prove the claim by induction on $r+\ell$. Let $c$ denote the class of 
$\Gamma$. By induction, we may reduce to the case where $C^{(c)}(\Gamma) = \ZZ g_r$.  

Since $\grisol \Gamma\otimes_\ZZ \QQ \cong \gr\Gamma\otimes_\ZZ \QQ$ is generated in degree $1$, 
there exist $s<t<r$ such that $\nu(t) = c-1$ and $[\bar{g}_s,\bar{g}_t] = n \bar{g}_r$, for 
some non-zero integer $n$. Then, $(g_s,g_t) = g_r^n$ and hence, by the previous calculations, 
one has $b_r(e_t,e_s) = -n$, while $b_r(e_s, e_t) = 0$. Therefore, by $(1)$,
 $[v_s,v_t] = n v_r$. 

{F}or any $k<r$, 
the image of $v_k$ in $\fraku/A v_r$ belongs to $C^{(\nu(k))}(\fraku/Av_r)$, 
by induction hypothesis.   
Thus, there exist a positive integer $m_k$ and $a_k\in A$ such that 
$$
m_k v_k - a_k v_r \in C^{\nu(k)}(\fraku). \leqno (4)
$$

Applying this to $k = t$ and using the fact that $v_r$ is central, one obtains that 
$$
m_t\,n \, v_r = [v_s, m_t v_t - a_t v_r] 
$$
belongs to $C^c(\fraku)$, and hence $v_r\in C^{(c)}(\fraku)$. In turn, this implies, 
by $(4)$, that 
$v_k \in C^{(\nu(k))}(\fraku)$, for each $k<r$. This proves the claim and  
 completes the proof of the proposition.

\subsection{Filtered Noetherian rings with the AR-property.}\label{AR} 

Let us recall several results about the homology of filtered Noetherian rings 
with the Artin-Rees property. Some basic references for this material are 
\cite{Sjo}, \cite{BD}, \cite{Gru}; see also \cite[Chap.\,I]{NVO} and 
\cite[\S\,1]{Do1}. (Note, however, that in \cite{Gru} 
the assertions in lines 8-12 of 2.8 and assertion (ii) of Theorem 3.3 are not correct; 
it is not difficult to provide counter-examples). 

Let $S$ be a left Noetherian ring. A sequence 
$\calI := \{I_1, I_2, \dots \}$ of two-sided ideals is said to be admissible if  
$I_1 \supseteq I_2 \supseteq \cdots$ and  
$I_j I_k \subseteq I_{j+k}$ for $j,k\geq 0$ (where one sets $I_0 = S$).   
Given such a sequence, let 
$$
\gr S := \bigoplus_{n\geq 0} I_n/I_{n+1} \quad\mbox{and}\quad 
\hatS := \mathop{\rm proj.lim.}_{n\geq 0} S/I_n
$$
be the associated graded ring and completion, respectively. 

Let $S$-filt denote the category of $\NN$-filtered 
left $S$-modules: objects are left 
$S$-modules $M$ equipped with a decreasing filtration 
$M = F^0 M\supseteq F^1 M\supseteq \cdots$ such that $I_n F^k M \subseteq F^{n+k}M$, 
and a morphism $f : M\to N$ between two such objects 
is an $S$-morphism which preserves the filtrations. Then $f$ induces a morphism 
of $\gr S$-modules $\gr f : \gr M \to \gr N$ and this defines a functor 
$\gr$ from $S$-filt to the category of $\NN$-graded $\gr S$-modules. 
Further, $f$ is called strict if one has $f(M)\cap F^k N = f(F^k M)$ for any $k$. 

An object $M$ of $S$-filt is called {\it separated\,} if\, $\bigcap_{n\geq 0} F^n M = \{0\}$, 
and {\it discrete\,} if $F^n M = \{0\}$ for some $n \geq 0$. 

The category $S$-filt is equipped with shift functors $s^n$, for $n\geq 0$, defined 
as follows. If $M$ is an object of $S$-filt, $s^n M = M$ as $S$-module but 
$F^p(s^n M) = F^{p-n}M$ for $p\geq 0$, with the convention that $F^k M = M$ if $k<0$. 
If $M$ is an $\NN$-graded $S$-module, the shifted module $s^n M$ is defined in 
an analogous manner. 

An object $L$ of $S$-filt is called filt-free 
if it a direct sum of shifted modules $s^{d(\lambda)}S$, for $\lambda$ running in some index set 
$\Lambda$. 
Then, $\gr L \cong \bigoplus_{\lambda\in\Lambda} s^{d(\lambda)} \gr S$.  

Let $M$ be an object of $S$-filt. Then a strict filt-free resolution is an  
$S$-module resolution 
$$
\cdots \longrightarrow L_1 \stackrel{f_1}{\longrightarrow} L_0 
\stackrel{f_0}{\longrightarrow} M \longrightarrow 0 \leqno (\calE)
$$
such that every $L_n$ is filt-free and every $f_n$ is a strict morphism in $S$-filt. 
By \cite[Lemmas 1,2]{Sjo}, the associated graded complex 
$(\gr \calE)$ is then a free $\gr S$-resolution of $\gr M$ and, conversely, 
if $S$ is complete with respect to $\calI$, 
any free $\gr S$-resolution of $\gr M$ can be obtained in this manner. 

\smallskip
Let us consider also the category filt-S of $\NN$-filtered {\it right\,} $S$-modules. 
All notions introduced previously for $S$-filt have, of course, their right-handed analogues. 
Now, if $N$ (resp. $M$) is an object of filt-$S$ (resp. $S$-filt), 
the abelian group $N\otimes_S M$ has a natural $\NN$-filtration, defined by
$$
F^n(N\otimes_S M) := 
\Ima\left(\sum\nolimits_{p+q=n} F^p N \otimes_S F^q M \to N\otimes_S M\right) .
$$
Moreover, it is easily seen that if either of $N$ or $M$ is a filt-free object, 
then the natural map 
$\gr N \otimes_{\gr S} \gr M \to \gr (N\otimes_S M)$ is an isomorphism. 

Therefore, if one considers a strict filt-free resolution $L_\bullet$ of, say, $M$, 
the filtration on $N\otimes_S L_\bullet$ induces a natural spectral sequence with 
$E_1$-term (in cohomological notation) 
$$
E_1^{p,-q} = H^{p-q}(\gr N\otimes_S \gr L_\bullet)_p = 
\Tor_{q-p}^{\gr S}(\gr N, \gr M)_p.
$$
Moreover, certain finiteness conditions ensure that 
this spectral sequence converges finitely to 
$\Tor_*^S(N,M)$. Firstly, 
by \cite[Lemma 2.(g)]{Sjo} or \cite[Th.\,2.9]{Gru}, 
one has the following 

\begin{propC} Assume that $S$ is complete with respect to the 
filtration $\calI$ and that $\gr S$ is left Noetherian. 
Let $M,N$ be objects of $S$-filt and filt-$S$, respectively, 
such that $M$ is separated and $\gr M$ finitely generated over $\gr S$,  
while $N$ is discrete. 
Then the spectral sequence above 
converges finitely to $\Tor_*^S(N,M)$. 
\end{propC}

\pf 
By the references cited above, any resolution of $\gr M$ by free 
$\gr S$-modules can be lifted to a strict filt-free resolution of $M$. 
Since $\gr M$ finitely generated over $\gr S$, which is left Noetherian,   
one deduces that $M$ admits a strict filt-free resolution 
$L_\bullet \to M \to 0$ such that each $L_n$ is finitely generated. 
As $N$ is assumed to be discrete, the filtration on 
$N\otimes_S L_\bullet$ is then discrete (and exhaustive) in each degree, and the 
proposition follows.

\medskip Secondly, the assumption that $S$ be complete can be relaxed if 
one assumes that 
the sequence $\calI = \{I = I_1\supseteq I_2 \supseteq \cdots\}$ has 
the left Artin-Rees property, {\it i.e.}, that $\calI$ satisfies the 
following : for any finitely generated left $S$-module $M$, any submodule $N \subseteq M$ 
and any $n\geq 0$, there exists 
$n'\geq n_0$ such that $N \cap I_{n'}M \subseteq I_n N$.

{F}or any left $S$-module $M$, let us denote by $\hatM$ its completion with respect to 
the filtration $\{I_n M\}$; it is an $\hatS$-module and there is a natural morphism 
of $\hatS$-modules $\tau_M : \hatS\otimes_S M \to \hatM$. 
As observed in \cite[Prop.\,3]{BD}, one has the following proposition, 
which is proved exactly as in the commutative $I$-adic case (see \cite[Chap.\,10]{AM}). 

\begin{propAR}  Assume that $S$ is left Noetherian and that $\calI$ satisfies the left 
AR-property. Then, $\tau_M$ is an isomorphism for any finitely generated left $S$-module 
$M$ and, therefore, \\
{\rm a)} $\hatS$ is flat as right $S$-module, \\
{\rm b)} for each $n$, $\hatS I_n = \Ker(\hatS \to S/I_n)$ is a two-sided ideal 
and hence $\{\hatS I_n\}$ is an admissible sequence in $\hatS$, \\
{\rm c)} the associated graded $\gr \hatS$ is isomorphic to $\gr S$. 
\end{propAR}
 
Thus, in particular, if $P_\bullet \to S/I \to 0$ is a resolution of $S/I$ by 
free $S$-modules, then $\hatS\otimes_S P_\bullet$ is a free 
$\hatS$-resolution of 
$$\hatS\otimes_S (S/I) = \widehat{S/I} = S/I.$$ 
Thus, for any right $\hatS$-module $N$, there is a natural isomorphism 
$$
\Tor_\bullet^\hatS(N, S/I) \cong \Tor_\bullet^S(N, S/I).
$$
This is the case, in particular, if $N$ is a right $S$-module with a discrete 
filtration. 
Therefore, one obtains the following theorem, which  
is essentially contained in \cite[Th.\,3.3$'$.(i)]{Gru}.

\begin{numth}\label{thAR} 
Let $S$ be a left Noetherian ring, $\calI$ an admissible sequence 
of ideals. Suppose that $\calI$ satisfies the left AR property and that 
$\gr S$ is left Noetherian. 
Let $N$ be a right $S$-module with a discrete filtration. 
Then there is a finitely convergent spectral sequence 
$$
E_1^{p,-q} = \Tor_{q-p}^{\gr S}(\gr N, S/I)_p \Rightarrow 
\Tor_{q-p}^\hatS(N, S/I) \cong \Tor_{q-p}^S(N, S/I).
$$
\end{numth}

For future use, let us derive the following equivariant version of the theorem. 
Let 
$\Lambda$ be a group of automorphisms of $S$ preserving the sequence $\calI$. 
Let $S\Lambda$ denote the smash product $S\# \ZZ\Lambda$, that is, 
$S\Lambda = S\otimes_\ZZ \ZZ\Lambda$ as $(S,\ZZ\Lambda)$-bimodule, the 
multiplication being defined by 
$$
(s\otimes \lambda)(s'\otimes \lambda') = s\lambda(s') \otimes \lambda\lambda'.
$$
Similarly, denote by $\hatS \Lambda$ the smash product $\hatS \# \ZZ\Lambda$. 
Observe that an $S\Lambda$-module is the same thing as an $S$-module $M$ 
equipped with an action of $\Lambda$ such that 
$\lambda sm = \lambda(s) \lambda m$, for $m\in M$, $s \in S$, $\lambda\in \Lambda$. 

For every $n\geq 0$, let $I'_n$ (resp. $\hatI'_n$) denote the left ideal of 
$S\Lambda$ (resp. $\hatS \Lambda$) generated by $I_n$; they are two-sided ideals and 
form an admissible sequence of $S\Lambda$ (resp. $\hatS \Lambda$). In both cases, 
the associated graded is isomorphic to $(\gr S)\Lambda := (\gr S)\# \Lambda$. 

\begin{numth}\label{thARequiv} 
With notation as above, let $N$ be a discrete object of 
$S\Lambda$-filt. There is a finitely convergent 
spectral sequence of $\Lambda$-modules 
$$
E_1^{p,-q} = \Tor_{q-p}^{\gr S}(\gr N, S/I)_p \Rightarrow \Tor_{q-p}^S(N, S/I).
$$
\end{numth} 

\pf First, $I' := (S\Lambda)I$ is a two-sided ideal of $S\Lambda$, and 
$S\Lambda\otimes_S (S/I) \cong S\Lambda/I'$. Then, 
by standard arguments, it suffices to prove that: {\it i\,})  $\hatS \Lambda$ is flat 
as right $S\Lambda$-module, and: {\it ii\,})  
$\hatS \Lambda\otimes_S (S/I) \cong S\Lambda/I'$. 

\vskip2pt 
But $\hatS \Lambda$ is isomorphic to $\hatS\otimes_S S\Lambda$ as 
$(\hatS, S\Lambda)$-bimodule, and to $S\Lambda \otimes_S \hatS$ as 
$(S\Lambda, \hatS)$-bimodule. This implies {\it i\,}) and {\it ii\,}).

\subsection{}\label{twospecseq} 

Let us return to the finitely generated, torsion free, nilpotent group $\Gamma$ 
and the associated unipotent algebraic group $U_A$. 
Recall the notation of subsections \ref{N-series}--\ref{same-gr}. 

It is known that $\ZZ\Gamma$ and $U_A(\fraku)$ are left and right Noetherian and 
have the left and right AR-property with respect to the filtration by the 
powers of the augmentation ideal, see, for example,  
\cite[Th.\,2.7 \& \S\,11.2]{Pas}, \cite{NG} and \cite[Th.\,1]{BD}.

Further, by \cite[Cor.\,3.5]{Hart}, one has $I^{(cn)}\subseteq I^{n}$, where $c$ is the class 
of $\Gamma$ (and also the class of $\fraku$), and a similar argument, using 
Proposition \ref{Vergne}.a) shows that 
$J^{(cn)} \subseteq J^{n}$. From this one deduces easily that the sequences 
$\{I^{(n)}\}$ and $\{J^{(n)}\}$ also have the left and right AR-property. 
In the sequel, we equip $\ZZ\Gamma$ and $U_A(\fraku)$ with these sequences, 
which we call $\calI$ and $\calJ$ respectively. 
By Theorem \ref{Hartley} and Proposition \ref{Vergne}, the associated graded rings are 
left and right Noetherian.  

\medskip 
Let $V$ be an $U_A$-module. Then $V$ is in a natural manner a representation of the Lie algebra 
$\fraku$ and of the abstract group $\Gamma$.   
Let $\calF$ be a finite sequence 
$V = F^0 V\supset \cdots \supset F^{s+1} V = \{0\}$ of $U_A$-submodules.  
Let us say that $\calF$ is an admissible filtration of $V$ if it is 
an $\calI$ (resp. $\calJ$) filtration of $V$ regarded as $\ZZ\Gamma$ (resp. $U_A(\fraku)$) 
module, {\it i.e.}, if    
for any $i,n\geq 0$, both $I^{(n)} (F^i V)$ and $J^{(n)} (F^i V)$ are contained in 
$F^{i+n} V$. 

\begin{lem} If $V$ is an $U_A$-module which is finite free over $A$, it admits 
an admissible filtration.
\end{lem}

\pf By the theorem of Lie-Kolchin applied to $V_\QQ$, one obtains 
that $V^U$, the submodule of invariants, is non-zero. 
Since 
$$
V^U = \{x\in V\mid \Delta_V(x) = x\otimes \varepsilon\},
$$
where $\Delta_V$ is the coaction defining the  comodule structure and $\varepsilon$ 
is the augmentation of $A[U]$, and since $V\otimes_A A[U]$ is a free $A$-module, 
one sees that $V/V^U$ is torsion-free, hence a free $A$-module. 

Therefore, if one sets $F_0 V = 0$ and defines inductively  $F_k V$ as the inverse 
image in $V$ of the $U$-invariants in $V/F_{k-1}V$, the sequence $\{F_k V\}$ 
is increasing strictly, as long as $F_k V \not= V$, and each $V/ F_k V$, if non-zero, 
 is a finite free $A$-module. Since $V$ is a Noetherian 
$A$-module, $F_N V = V$ for some $N$. Setting $F^i V = F_{N-i} V$, 
it is easily seen that, for any $i,n\geq 0$, both $I^n (F^i V)$ and $J^n (F^i V)$ 
are contained in $F^{i+n} V$. Further, since every $F^i V/F^{i+n} V$ is torsion-free, 
 one obtains that $\{F^i V\}$ is an admissible filtration of $V$.

\smallskip 
Then, one deduces from the results of \ref{AR} the following theorem. There are, obviously, 
equivariant versions; we leave their formulation to the reader.   

\begin{th} Let $V$ be an $U_A$-module which is finite free over $A$ and let  
$\calF$ be any admissible filtration on $V$. 
Then there are two finitely convergent spectral sequences:  
$$
\displaylines{ 
\phantom{i}{\it i\,}) \qquad\qquad
E_1^{p,-q} = H_{q-p}(\grisol \Gamma, \gr_\calF V)_p \Rightarrow  
H_{q-p}(\Gamma, V), 
\cr
 {\it ii\,}) \qquad\qquad
E_1^{p,-q} = H_{q-p}(\grisol \fraku, \gr_\calF V)_p \Rightarrow 
 H_{q-p}(\fraku, V). 
}
$$    
\end{th}

\subsection{}\label{f_I}
Finally, let us return to the setting of Sections 1 and 2. 
The unipotent group $U_P^-$ is defined over $\ZZ$.   
Let $\Gamma := U_P^-(\ZZ)$; it is, clearly, a torsion-free nilpotent group.   
 
Let $\{H_\alpha\}_{\alpha\in\Delta} \cup \{X_\beta\}_{\beta\in R}$ be a Chevalley basis of 
$\frakg$. For each $\beta\in R$, let $U_\beta$ be the corresponding root subgroup and 
let $\theta_\beta$ be the isomorphism $\GG_a\to U_\beta$ such that $d\theta_\beta(1) = X_\beta$. 
Set $I := \Delta\setminus R_L^+$ and let $f_I : \ZZ R \to \ZZ$ be the additive function 
which coincides on the basis $\Delta$ with the negative of the characteristic function of $I$. 
That is,  
$$
f_I(\alpha) = \left\{\begin{array}{ll} -1 & \mbox{if $\alpha\in I$};\\
0 & \mbox{if $\alpha\in \Delta\cap R_L^+$}.\end{array}\right.
$$
Choose a numbering $\alpha_1,\dots,\alpha_r$ of the elements of $R^+\setminus R_L^+$ 
such that $f_I(\alpha_i)\leq f_I(\alpha_j)$ if $i\leq j$. The 
multiplication map induces an isomorphism of $\ZZ$-schemes 
$$
U_{\alpha_1}\times\cdots\times U_{\alpha_r} \stackrel{\cong}{\to} U_P^- .
$$
Moreover, it follows from the commutation formulas in \cite[Lemma 15]{St} 
or \cite[3.2.3--3.2.5]{BT} that, for any $s=1,\dots,r$, $U_{\alpha_s}\cdots U_{\alpha_r}$ 
is a closed, normal subgroup of $U_P^-$. One deduces that the $g_i := \theta_{\alpha_i}(1)$ 
form a system of canonical parameters of $\Gamma$, that $U_P^-$ is the algebraic group 
associated in \ref{Hall} to $\Gamma$, and that the basis $\{v_1,\dots,v_r\}$ of 
$\fraku_P^-$ identifies with $\{X_{\alpha_1},\dots,X_{\alpha_r}\}$.

\begin{lem} One has $\fraku_P^- \cong \grisol\fraku_P^-$.
\end{lem} 

\pf Since $T$ acts on $\fraku_P^-$ by Lie algebra automorphisms, $\fraku_P^-$ has 
a structure of graded Lie algebra given by the function $f_I$. 
That is, if one sets, for $i\geq 1$, 
$$
\fraku_P^-(i) := \bigoplus_{{\alpha\in R^- \atop f_I(\alpha) =i}} \frakg_\alpha,
$$
then 
$$
\fraku_P^- = \bigoplus_{i\geq 1} \fraku_P^-(i)\quad\mbox{and}\quad [\fraku_P^-(i),\fraku_P^-(j)] 
\subseteq \fraku_P^-(i+j).
$$
 Therefore, the lemma will follow if we show that 
$C^{(i)}(\fraku_P^-) = \fraku_P^-(\geq i)$, where $\fraku_P^-(\geq i)$ is defined in the obvious 
manner. Clearly, $C^i(\fraku_P^-) \subseteq \fraku_P^-(\geq i)$ and, since 
$\fraku_P^-/\fraku_P^-(\geq i)$ is torsion-free, one obtains that 
$C^{(i)}(\fraku_P^-) \subseteq \fraku_P^-(\geq i)$. 

\smallskip
In order to prove the converse inclusion, it suffices to prove that 
$\fraku_{P,\QQ}^-(i)\subseteq C^i(\fraku_{P,\QQ}^-)$, where $\fraku_{P,\QQ}^- = 
\fraku_P^-\otimes_\ZZ \QQ$. But this follows from Kostant's theorem about the homology 
of $\fraku_{P,\QQ}^-$ (\cite[Cor.\,8.1]{Ko1}). Indeed, $\fraku_{P,\QQ}^-$ is generated by any 
subspace supplementary to $[\fraku_{P,\QQ}^-,\fraku_{P,\QQ}^-]$. But, by the cited result 
of Kostant, one has 
$$
\fraku_{P,\QQ}^- / [\fraku_{P,\QQ}^-,\fraku_{P,\QQ}^-] = H_1(\fraku_{P,\QQ}^-,\QQ) 
\cong \bigoplus_{\alpha\in I} V_\QQ^L(-\alpha),
$$
and the R.H.S. identifies with $\fraku_{P,\QQ}^-(1)$. Therefore, 
$\fraku_{P,\QQ}^-$ is generated by $\fraku_{P,\QQ}^-(1)$. Then, by induction on $i$, 
one obtains easily that $\fraku_{P,\QQ}^-(i)\subseteq C^i(\fraku_{P,\QQ}^-)$ for any $i$. 
The lemma is proved. 

\medskip Recall the integers $\nu(i)$ introduced in the proof of Proposition \ref{same-gr}. 
From this proposition and the previous lemma (and their proofs), one deduces the following

\begin{cor} There is an isomorphism of graded Hopf algebras 
$\grisol\ZZ\Gamma \cong U(\fraku_P^-)$, 
under which each $\overline{g_i-1}$ correponds to $X_{\alpha_i}$. 
Further, for  $i=1,\dots,r$, one has $\nu(i) = f_I(\alpha_i)$.
\end{cor}

\subsection{}\label{grV}  
For any $\lambda\in X^+$, set 
$$
V_\ZZ(\lambda)(i) := \bigoplus_{{\mu\in X \atop f_I(\mu-\lambda) = i}} V_\ZZ(\lambda)_\mu,
$$
where the subscript $\mu$ denotes the $\mu$-weight space.
Then, each $V_\ZZ(\lambda)(i)$ is an $L$-submodule and there is an isomorphism 
of $L$-modules 
$$
V_\ZZ(\lambda) \cong \bigoplus_{i\geq 0} V_\ZZ(\lambda)(i).
$$
Set $F^k V_\ZZ(\lambda) := \bigoplus_{i\geq k} V_\ZZ(\lambda)(i)$; this defines a filtration 
$\calF$ of $V_\ZZ(\lambda)$ by $P^-$-submodules, such that the associated graded 
is isomorphic to $V_\ZZ(\lambda)$ as $L$-module. 

\begin{prop} One has $I^{(n)} F^k V_\ZZ(\lambda)\subseteq F^{n+k} V_\ZZ(\lambda)$, 
and $\gr_\calF V_\ZZ(\lambda) \cong V_\ZZ(\lambda)$ as   
representations of $\grisol\ZZ\Gamma\cong U_\ZZ(\fraku_P^-)$. 
\end{prop} 

\pf 
For $i = 1,\dots,r$ and $n\geq 0$, set 
$$
u_i^{(n)} := g_i^{-[(n+1)/2]}\, (g_i-1)^n , 
$$
where $[x]$ denotes the greatest integer not greater than $x$,  
and observe that $u_i^{(n)} \equiv (g_i-1)^n$ modulo $I^n$. 
Further, for ${\bf j}\in \NN^r$, set 
$$
u({\bf j}) := u_1^{(j_1)} \cdots u_r^{(j_r)} \quad \mbox{and}\quad 
\nu({\bf j}) = \sum_i j_i \nu(i). 
$$
Then, by \cite{Hart}, Theorem 3.2\,(i) and Lemma 3.1, the elements $u({\bf j})$ satisfying 
$\nu({\bf j})\geq n$ form a $\ZZ$-basis of $I^{(n)}$, for every $n\geq 0$. 

From this one deduces that, in order to prove the proposition, it suffices to prove that, 
for any $v\in F^k V_\ZZ(\lambda)$ and $i=1,\dots,r$, one has 
$$
(g_i-1) v - X_{\alpha_i} v \in F^{k+\nu(i)+1} V_\ZZ(\lambda). \leqno (*) 
$$

The distribution algebra $\dist(U_P^-)$ has a $\ZZ$-basis formed by 
the ordered products 
$$
X_{\alpha_1}^{(m_1)}\cdots X_{\alpha_r}^{(m_r)},\quad\mbox{for }(m_1,\dots,m_r)\in\NN^r,
$$
where the elements $X_\beta^{(m)}$ satisfy $X_\beta^m = m!\, X_\beta^{(m)}$ for every 
$m\geq 0$. Further, the structure of $\ZZ[G]$-comodule on $V_\ZZ(\lambda)$ is such that, 
for any ring $\Omega$, any $t\in\Omega$ and $v\in V_\Omega(\lambda)$, and any 
root $\alpha$, one has 
$$
\theta_\alpha(t)v = \sum_{m\geq 0} t^m X_\alpha^{(m)} v,
$$
where the R.H.S. is in fact a finite sum. Since $g_i = \theta_{\alpha_i}(1)$ and 
since each $X_{\alpha_i}^{(m)}$ has weight $m\alpha_i$ for the adjoint action of $T$, 
this immediately implies formula $(*)$. The proposition is proved.

\subsection{}\label{mainthree} 
We can now prove Theorem \ref{introPickel} of the Introduction. 
The discrete group $\Lambda = L(\ZZ)$ normalizes $\Gamma = U_P^-$ and, hence, 
preserves the isolated powers of the augmentation ideal of $\ZZ\Gamma$. 
Therefore, by the equivariant version of Theorem \ref{twospecseq}\,{\it i\,}), 
combined with Proposition \ref{grV}, there is a 
finitely convergent spectral sequence of $L(\ZZ)$-modules 
$$
H_*(\fraku_P^-, V_\ZZ(\lambda)) \cong H_*(\grisol\ZZ\Gamma, V_\ZZ(\lambda)) 
\Rightarrow H_*(\Gamma, V_\ZZ(\lambda)). \leqno (1)
$$
It is, clearly, compatible with flat base change. Thus, for any prime integer $p$, 
one has a finitely convergent spectral sequence
$$
H_*(\fraku_P^-, V_\Zp(\lambda)) \cong H_*(\grisol\ZZ\Gamma, V_\Zp(\lambda))  
\Rightarrow H_*(\Gamma, V_\Zp(\lambda)). \leqno (2)
$$
Moreover, it is not difficult to check, by standard arguments, that the 
natural structure of $L(\ZZ)$-module on 
$H_*(\grisol\ZZ\Gamma, V_\Zp(\lambda))$ considered in Theorem \ref{thARequiv} 
is the restriction to $L(\ZZ)$ of the natural structure of $L$-module on 
$H_*(\fraku_P^-, V_\Zp(\lambda))$. Therefore, if $\lambda$ is $p$-small then, 
by Theorem \ref{kostant}, one obtains an isomorphism of 
$L(\ZZ)$-modules 
$$
H_i(\grisol\ZZ\Gamma, V_\Zp(\lambda)) \cong H_i(\fraku_P^-, V_\Zp(\lambda)) \cong 
\bigoplus_{w\in W^L(i)} V_\Zp^L(w\cdot\lambda),
$$
for every $i\geq 0$. In particular, $H_*(\grisol\ZZ\Gamma, V_\Zp(\lambda))$ is a 
free $\Zp$-module. 

Finally, it is well-known that $\fraku_P^-\otimes \QQ$ is isomorphic to the 
Malcev-Jennings Lie algebra of $\Gamma$; this follows, for example, from the 
proof of \cite[Lemma 1.9]{LP}. Therefore, by a result of Pickel \cite[Th.\,10]{Pic}, 
there is an isomorphism of graded vector spaces 
$$
H_\bullet(\fraku_P^-, V_\QQ(\lambda)) \cong H_\bullet(\Gamma, V_\QQ(\lambda)).
$$
This implies that the abutment of the spectral sequence in $(2)$ has the same rank 
over $\Zp$ as the $E_1$-term. Since the latter is a free $\Zp$-module, one deduces that 
the spectral sequence degenerates at $E_1$. Therefore, we have obtained the 
following 

\begin{th} Let $\lambda\in X^+\cap \barCp$. Then, for each $n\geq 0$,  
$H_n(U_P^-(\ZZ), V_\Zp(\lambda))$ has a finite, natural $L(\ZZ)$-module filtration    
such that 
$$
gr H_n(U_P^-(\ZZ), V_\Zp(\lambda)) \cong \bigoplus_{w\in W^L(n)} V_\Zp^L(w\cdot\lambda).
$$
\end{th} 

By the universal coefficient theorem, one then obtains a similar result over $\Fp$. 
Finally, by an argument similar to the one in \ref{cohom-u}, one obtains the

\begin{cor} Let $\lambda\in X^+\cap \barCp$. Then, for each $n\geq 0$,  
$H^n(U_P(\ZZ), V_\Fp(\lambda))$ has a finite, natural $L(\ZZ)$-module filtration   
such that 
$$
gr H^n(U_P(\ZZ), V_\Fp(\lambda)) \cong \bigoplus_{w\in W^L(n)} V_\Fp^L(w\cdot\lambda).
$$
\end{cor}

\subsection{}\label{p-Lie} 

Let us derive in this subsection a corollary about the $p$-Lie algebra associated with 
the $p$-lower central series of $\Gamma$. (This result will not be used in the sequel). 

Let $\calF$ be a decreasing sequence 
$\Gamma = F^1\Gamma \supseteq F^2\Gamma \supseteq \cdots$ of normal subgroups 
of $\Gamma$.  
It is called an $N_p$-sequence if it is an $N$-sequence and $x\in F^i\Gamma$ 
implies that $x^p \in F^{pi}\Gamma$. 
In this case, $\gr_\calF \Gamma$ is a graded $p$-Lie algebra, 
see \cite[Chap.\,I, Cor.\,6.8]{Laz}  or \cite[Chap.\,II, \S\,5, Ex.\,10]{Bki}.

For our purposes, it is convenient to define the $p$-lower central series 
$\{F_p^n\Gamma\}_{n\geq 1}$ as follows. 
Denoting by $I_\Fp$ the augmentation ideal of $\Fp\Gamma$, set 
$$
F_p^n \Gamma := \{x\in \Gamma \mid  x-1\in I_\Fp^n\}.
$$
This is an $N_p$-sequence (see \cite[Lemma 3.3.1]{Pas}), and we denote 
the associated graded $p$-Lie algebra by  
$\gr_p^\bullet \Gamma$. 

The $n$-th term $F_p^n\Gamma$ of the $p$-lower central series is sometimes defined as the 
subgroup of $\Gamma$ generated by all elements $x^{p^s}$ satisfying $p^s\omega(x)\geq n$, 
where $\omega(x)$ denotes the 
largest integer $i$ such that $x\in C^i(\Gamma)$. That the two definitions agree is due to Lazard 
\cite[Chap.\,I, Th.\,5.6 \& 6.10]{Laz} and Quillen \cite{Qui}, see also \cite[\S\,11.1]{Pas}. 

\smallskip 
Let us denote by $\calLie_\Fp$ the category of Lie algebras over $\Fp$, by $\pcalLie_\Fp$ 
the subcategory of $p$-Lie algebras, and by $\grcalLie_\Fp$ and 
$\pgrcalLie_\Fp$, 
respectively, the subcategories of graded and graded $p$-Lie algebras over $\Fp$. 
The forgetful functor $\pcalLie_\Fp \to \calLie_\Fp$ has a left adjoint, denoted by 
$\pcalL$; it takes $\grcalLie_\Fp$ to $\pgrcalLie_\Fp$. 

\begin{cor} Let $\Gamma$ be a finitely generated, torsion-free, nilpotent group, say of class 
$c$. Suppose that $\oplus_{i=1}^c\, C^{(i)}(\Gamma)/C^i(\Gamma)$ has no $p$-torsion. 
Then, there is 
an isomorphism of graded $p$-restricted Lie algebras 
$$
\gr_p^\bullet\,\Gamma \cong \pcalL(\gr\Gamma\otimes \Fp).
$$
\end{cor}

\pf The hypothesis implies easily that $\gr\Gamma\otimes\Fp \cong\grisol\Gamma\otimes\Fp$. 
Moreover, it follows from the proof of \cite[Th.\,3.2\,(i)]{Hart} that every 
$I^{(n)}/I^n$ has no $p$-torsion. This implies that, inside $\Fp\Gamma$, one has the identifications 
$I^{(n)}\otimes\Fp = I^n\otimes\Fp = I_\Fp^n$. 
One deduces from this, coupled with Theorem \ref{Hartley},  
the isomorphisms 
$$
\displaylines{ 
\gr\Fp\Gamma\cong (\grisol\ZZ\Gamma)\otimes\Fp \cong U_\ZZ(\grisol\Gamma)\otimes\Fp 
\hfill\cr\hfill
\cong U_\Fp(\grisol\Gamma\otimes\Fp) \cong U_\Fp(\gr\Gamma\otimes\Fp).
}
$$
On the other hand, by Quillen \cite{Qui}, 
$\gr\Fp\Gamma$ is isomorphic as graded Hopf algebra to $U_\Fp^{res}(\gr_p^\bullet\Gamma)$, 
the restricted enveloping algebra of the $p$-Lie algebra $\gr_p^\bullet\Gamma$. 

Recall that $U_\Fp^{res}$, the restricted enveloping algebra functor, is left adjoint to 
the forgetful functor $\calAs_\Fp\to \pcalLie_\Fp$, where $\calAs_\Fp$ denotes the 
category of associative $\Fp$-algebras (with unit), while the usual enveloping algebra 
functor is left adjoint to the forgetful functor $\calAs_\Fp\to \calLie_\Fp$. Thus, 
since the adjoint of a composite is the composite of the adjoints, one has 
$U_\Fp(L) \cong U_\Fp^{res}(\pcalL(L))$, for any $\Fp$-Lie algebra $L$. 

Therefore, one obtains an isomorphism of graded Hopf algebras 
$$
U_\Fp^{res}(\pcalL(\gr\Gamma\otimes\Fp)) \cong U_\Fp^{res}(\gr_p^\bullet\Gamma).
$$ 
Taking primitive elements, this gives, 
by the theorem of Milnor-Moore \cite[Th.\,6.11]{MM},   
an isomorphism of graded $p$-Lie algebras $\pcalL(\gr\Gamma\otimes\Fp) \cong 
\gr_p^\bullet\Gamma$. 
The corollary is proved. 

\medskip\noindent {\it Remark}. It is easy to see that the torsion primes in 
$\oplus_{i=1}^c C^{(i)}(\Gamma)/C^i(\Gamma)$ and in $\gr \Gamma$ are the same. 
Presumably, it should not be difficult to extract from the proof 
of Proposition \ref{same-gr} that the torsion primes in $\gr\fraku$ are also the same.

\section{Standard and BGG complexes for distribution algebras}

\subsection{}
As in subsection \ref{standard}, there is defined a complex 
$$
\cdots \to \calU(G)\otimes_{\calU(P)} \Lambda^2(\frakg/ \frakp) 
 \stackrel{d_2^\frakp}{\to} \calU(G)\otimes_{\calU(P)} (\frakg/ \frakp)   
\stackrel{d_1^\frakp}{\to} \calU(G)\otimes_{\calU(P)} \ZZ 
\stackrel{\varepsilon}{\to} \ZZ\to 0,  
$$
the differentials being defined by the same formula as in \ref{standard}. Note, however, 
that this complex is {\it not\,} exact. We shall denote it by 
$\calS_\bullet(G,P)$. 

More generally, let $V$ be a $G$-module and let $V_\vert{}_P$ denote $V$ regarded as 
an $\calU(P)$-module. Then one obtains, as in \ref{standard}, a complex of 
$\calU(G)$-modules 
$$
\displaylines{ 
\cdots \to \calU(G)\otimes_{\calU(P)} (\Lambda^2(\frakg/ \frakp) \otimes V\vert{}_P) 
 \stackrel{d_2}{\to} \calU(G)\otimes_{\calU(P)} (\frakg/ \frakp \otimes V\vert{}_P)  
\hfill\cr\hfill  
\stackrel{d_1}{\to} \calU(G)\otimes_{\calU(P)} V\vert{}_P 
\stackrel{\varepsilon}{\to} V\to 0. 
}  
$$
We shall call it the standard complex of $V$ relative to the pair $(\calU(G), \calU(P)$,  
and denote it by $\calS_\bullet(G,P,V)$. When $V = V_\ZZ(\lambda)$, we shall denote it 
simply by $\calS_\bullet(G,P,\lambda)$. 

\smallskip 
Further, as in \ref{defVerma}, let us define, for any $\xi\in X_L^+$, the 
generalized Verma module (for $\calU(G)$ and $\calU(P)$) 
$$
\calM_P(\xi) := \calU(G)\otimes_{\calU(P)} V_\ZZ^L(\xi).
$$
Moreover, for any commutative ring $A$, set $\calM_P^A(\xi) := \calM_P(\xi)\otimes A$ 
and, for any $\lambda \in X^+$, let 
$$
\calS_\bullet^A(G,P,\lambda) := \calS_\bullet(G,P,\lambda)\otimes A. 
$$

\subsection{}\label{th-BGGcomplex}
Our aim in this section is to prove the following theorem. 

\begin{th} Suppose that $\fraku_P^-$ is abelian. Let $\lambda\in X^+\cap\barCp$. Then 
the standard complex $\calS^\Zp_\bullet(G,P,\lambda)$ contains as a direct summand 
a subcomplex $\calC^\Zp_\bullet(G,P,\lambda)$ such that, for $i\geq 0$, 
$$
\calC_i^\Zp(G,P,\lambda) 
\cong 
\bigoplus_{w\in W^L(i)} \calM_P^\Zp(w\cdot\lambda).
$$
\end{th}

\subsection{The case $\lambda = 0$.}\label{lambda=0} 
As observed by Faltings-Chai in \cite[VI.5, p.230]{CF}, if $\fraku_P^-$ is abelian then 
$H_\bullet(\fraku_P^-)\cong \Lambda^\bullet(\fraku_P^-)$ and hence, 
by a result of Kostant \cite[\S\,8.2]{Ko1}, 
the composition  factors
of $\Lambda^i(\frakg/\frakp)_\QQ$ are exactly the $V_\QQ^L(w\cdot 0)$, for 
$w\in W^L(i)$, each occuring with multiplicity one. 

Thus, by Lemma \ref{donkin} and Corollary \ref{vanish}, applied to $L$, 
each $\Lambda^i(\frakg/\frakp)_\Zp$ is the direct sum of the Weyl modules 
$V_\Zp^L(w\cdot 0)$, for $w\in W^L(i)$. 
It follows that 
$$ 
\calS_i^\Zp(G,P) \cong \bigoplus_{ w\in W^L(i)} 
\calM_P^\Zp(w\cdot 0). \leqno (*)
$$ 
This proves the sought-for result when $\lambda = 0$ and $p \geq h-1$ ($h$ being the Coxeter 
number, see \ref{cox}).

\subsection{The general case.}

Now, let $\lambda\in X^+\cap\barCp$. 
First, since $\calS_\bullet^\Zp(G,P,\lambda) =
\calS_\bullet^\Zp(G,P)\otimes V(\lambda)$, it follows from 
\ref{lambda=0}\,$(*)$ and the tensor identity  
(\cite[Prop.\,1.7]{GL}) that, for $i\geq 0$, 
$$ 
\calS_i^\Zp(G,P,\lambda) \cong \bigoplus_{w\in W^L(i)}  
\calU_\Zp(G)\otimes_{\calU_\Zp(P)} \left( V_\Zp^L(w\rho-\rho) \otimes V_\Zp(\lambda) \right). 
\leqno (1)
$$
Let $\calS_w^\Zp(G,P,\lambda)$ denote the summand corresponding to $w$ in the R.H.S. 
Then 
$$
\calS_\bullet^\Zp(G,P,\lambda) = \bigoplus_{w\in W^L} 
\calS_w^\Zp(G,P,\lambda), \leqno (2)
$$
each $\calS_w^\Zp(G,P,\lambda)$ occuring in degree $\ell(w)$. 

Recall the notation $U'_A(\frakg)$ from \ref{Carter-Lusztig}. 
Since $U'_\Zp(\frakg)\subset \calU_\Zp(G)\subset U_\QQ(\frakg)$, one 
deduces that $U'_\Zp(\frakg)^G$ is contained in the center of $\calU_\Zp(G)$. 
Therefore, using exactly the same arguments as in \ref{blocks} and \ref{weakBGG}, 
one obtains that:  
a) $\calS_\bullet^\Zp(G,P,\lambda)$ contains as a direct summand 
the subcomplex 
$$
\calS_\bullet^\Zp(G,P,\lambda)_{\barchi_{\lambda,p}} \cong 
\bigoplus_{w\in W^L} 
\calS_w^\Zp(G,P,\lambda)_{\barchi_{\lambda,p}}, \leqno (3) 
$$
b) the latter has a filtration with associated graded 
$$
\gr \calS_\bullet^\Zp(G,P,\lambda)_{\barchi_{\lambda,p}} \cong 
\bigoplus_{y\in W^L} \calM_P^\Zp(y\cdot\lambda), \leqno (4)
$$
each $\calM_P^\Zp(y\cdot\lambda)$ occuring in degree $\ell(y)$, and c)  
for each $w\in W^L$, the sub\-quotients occuring in the filtration of 
$\calS_w^\Zp(G,P,\lambda)_{\barchi_{\lambda,p}}$ are those  
$\calM_P^\Zp(y\cdot\lambda)$ such that 
$V_\QQ^L(y\cdot\lambda)$ is a composition factor of the $L_\QQ$-module 
$V_\QQ^L(w\cdot 0)\otimes V_\QQ(\lambda)$. 

But, it is well-known that, necessarily, $y = w$. This may be deduced, for example, 
from \cite[Satz 2.25]{Ja-LN}. For the convenience of the reader, let us record a proof. 
First, it is well-known that any composition factor of the $L_\QQ$-module  
$V_\QQ^L(w\cdot 0)\otimes V_\QQ(\lambda)$ has the form $V_\QQ^L(w\cdot 0 +\nu)$, for 
some weight $\nu$ of $V_\QQ(\lambda)$, see, for example, 
\cite[\S\,24, Ex.\,12]{Hu} or, better, the proof of Cor.\,4.7 in \cite{An}. 
Secondly, for such a $\nu$, suppose that 
$w\cdot 0 + \nu = y\cdot\lambda$, for some $y\in W$. Then, 
$$
y^{-1}w\rho - \rho = \lambda - y^{-1}\nu.
$$
Let $\theta$ denote this weight. Since $y^{-1}w\rho$ (resp. $y^{-1}\nu$) is a weight of 
$V_\QQ(\rho)$ (resp. $V_\QQ(\lambda)$), one has $\theta\in -\NN R^+$ 
(resp. $\theta\in \NN R^+$) and, therefore, $\theta = 0$. Thus, since the stabilizer of 
$\rho$ in $W$ is trivial, $y = w$. 

Then, by comparing $(3)$ and $(4)$, one deduces that 
$\calS_w^\Zp(G,P,\lambda)_{\barchi_{\lambda,p}} 
\cong \calM_P^\Zp(w\cdot\lambda)$, for every $w\in W^L$. 
This completes the proof of Theorem \ref{th-BGGcomplex}.

\section{Dictionary}\label{dict} 

 Let $G=GSp(2g)_\ZZ$ be the split reductive Chevalley group over $\ZZ$ defined by 
${}^tXJX=\nu\cdot J$
where $J$ is given by $g\times g$-blocks  
$$
J=\left(
\begin{array}{cc}0_g&\begin{array}{cc}&\adots\\1&\end{array}\\
\begin{array}{cc}&\adots\\ -1&\end{array}&0_g\end{array}
\right)
$$
Let $B = TN$, resp. $Q = MU$, be the Levi decomposition of 
the upper triangular subgroup of $G$,  resp. of the Siegel parabolic, 
{\it i.e.}, the maximal parabolic associated to $\alpha$,  the longest simple
root for $(G,B,T)$, so $M=L_I$ where $I=\Delta\backslash\{\alpha\}$. 
Note that the unipotent radical of $Q$ is abelian. 

The group of characters $X$ of $T$ is identified to the sublattice 
$$
\{(a_g,\cdots,a_1;c)\in \ZZ^g\times\ZZ \mid c\equiv a_g+\ldots +a_1\,\,mod.2\}
$$
of $\ZZ^{g+1}$ in the following manner. The character $(a_g,\cdots,a_1;c)$ 
is defined by 
$$
diag (t_g,\ldots,t_1,x\cdot t_1^{-1},\ldots,x\cdot t_g^{-1})\mapsto
t_g^{a_g}\cdot\ldots\cdot t_1^{a_1}\cdot x^{(c-a_1-\ldots-a_g)/2}.
$$
The half-sum $\rho$  of the positive roots of $G$ is then $\rho=(g,\ldots,1;0)$.
If $(\varepsilon_g,\ldots,\varepsilon_1)$ is the canonical basis of $\ZZ^g$, the highest coroot 
$\gamma^\vee$ of $G$ is $\varepsilon_g+\varepsilon_{g-1}$. The condition
 $\langle\lambda+\rho,\gamma^\vee\rangle \leq p$,  
reads therefore:
$$a_g+a_{g-1}+g+(g-1) \leq p$$

 The lattice of weights $P(R)$ coincides with $X$. The cone $X^+\subset X$ of
dominant weights of $G$ is then given by the conditions $a_g\geq\ldots a_1\geq 0$.  For a character
$\phi=(a_g,\ldots,a_1;c)$ we define its degree as $\vert \phi\vert=\sum_{i=1}^g a_i$; the dual
$\hat{\phi}=(a_g,\ldots,a_1;-c)$ of $\phi$ has same degree
$\vert\hat{\phi}\vert=\vert\phi\vert$. 
 Note that $\vert\rho\vert=g(g+1)/2$. So,
$$\langle\lambda+\rho,\gamma^\vee\rangle\leq \vert\lambda+\rho\vert$$
with equality for $g\leq 2$.
 Let
${\bf V}=\langle e_g,\ldots,e_1,e_1^*,\ldots,e_g^*\rangle$ be the standard
$\ZZ$-lattice on which
$G$ acts; given two vectors $v,w\in \bV$, we write $<v,w>={}^tvJw$ for their symplectic product. 
$Q$ is the stabilizer of the standard lagrangian lattice $\bW=\langle e_g,\ldots,e_1 \rangle$; 
we have
${\bf V}=\bW\oplus
\bW^*$; $M=L_I$ is the stabilizer of the decomposition $(\bW,\bW^*)$;  one has
$M\cong GL(g)\times \GG_m$. Let
$B_M=B\cap M$ be the standard Borel of $M$.

\smallskip Recall from \ref{lattices} the definition of admissible lattices and, 
for $\lambda\in X^+$, the 
$\ZZ$-lattices $V(\lambda)_{min}$ and $V(\lambda)_{max}$. 

Let $\lambda\in X^+$ and $n= \vert\lambda \vert$; for any $(i,j)$ with $1\leq i<j\leq n$, let 
 $\phi_{i,j}:\bV^{\otimes n}\rightarrow\bV^{\otimes (n-2)}$
denote the contraction given by
$$
v_1\otimes\ldots\otimes v_n\mapsto\langle v_i,v_j\rangle v_1\otimes\ldots\otimes
\hat{v}_i\otimes\ldots\otimes\hat{v}_j\otimes\ldots\otimes v_n,
$$ 
and let $\bV^{<n>}$ 
be the submodule of
$\bV^{\otimes n}$ defined as intersection of the kernels of the $\phi_{i,j}$.
By applying the Young symmetrizer
$c_\lambda=a_\lambda\cdot b_\lambda$ (see
\cite{FH} 15.3 and 17.3) to  $\bV^{<n>}$, one obtains an admissible $\ZZ$-lattice  
$V(\lambda)_{Young}$ in $V_\QQ(\lambda)$. 

Then, by Corollary \ref{irred}, one has the 
 
\begin{cor} For any $p$-small weight $\lambda\in X^+$, one has canonically
$$
V(\lambda)_{min}\otimes {\Zp} = V(\lambda)_{Young} \otimes {\Zp}
= V(\lambda)_{max}\otimes {\Zp}. 
$$

Moreover, a similar result holds for a
 weight $\mu\in X^+_M$ of $M$, provided it is $p$-small for $M$.
\end{cor}

} 

\end{document}